\documentclass[ijoc,nonblindrev]{informs3} 

\OneAndAHalfSpacedXII 

\usepackage{natbib}
 \bibpunct[, ]{(}{)}{,}{a}{}{,}%
 \def\bibfont{\small}%
\TheoremsNumberedThrough     

\EquationsNumberedThrough    

\usepackage{algorithm}  
\usepackage{algorithmic} 

\usepackage{array}
\usepackage[dvipsnames]{xcolor}
\usepackage{tikz}
\usepackage{wrapfig}
\usepackage{microtype}
\usepackage{graphicx}
\usepackage{booktabs}
\usepackage{hyperref}

\usepackage{amssymb}
\usepackage{amsfonts}
\usepackage{amsmath}
\usepackage{latexsym}
\usepackage{epsfig}

\usepackage{multirow}
\usepackage{verbatim}
\usepackage{threeparttable}
\usepackage{color, colortbl}
\usepackage{multicol}
\usepackage{xr} 
\usepackage{caption}
\usepackage{subcaption}

\usepackage{todonotes}

\def\endproof{\hfill $\Box$ \vskip 0.4cm}

\newcommand{\expnumber}[2]{{#1}\mathrm{e}{#2}}
\usepackage{cleveref}

\usepackage{tikz}
\usetikzlibrary{shapes,positioning}
\tikzstyle{inout}=[trapezium, trapezium left angle=60, trapezium right angle=120, draw] 
\tikzstyle{end}=[rectangle, rounded corners, draw]   
\tikzstyle{endn}=[rounded rectangle, draw]   
\tikzstyle{exec}=[rectangle, draw]    
\tikzstyle{decide}=[kite, kite vertex angles=120, draw]   

\crefname{section}{Section}{Section}
\crefname{figure}{Figure}{Figure}
\crefname{table}{Table}{Table}
\crefname{equation}{}{}

\begin{document}


\RUNAUTHOR{Ge, Wang, Xiong, Ye}

\RUNTITLE{From an Interior Point to a Corner Point: Smart Crossover}

\TITLE{From an Interior Point to a Corner Point: Smart Crossover}

\ARTICLEAUTHORS{
\AUTHOR{Dongdong Ge}
\AFF{Antai School of Economics and Management, Shanghai Jiao Tong University, Shanghai, 200030, China. \EMAIL{ddge@sjtu.edu.cn}}
\AUTHOR{Chengwenjian Wang}
\AFF{Department of Industrial and Systems Engineering, Univesity of Minnesota, Minneapolis, MN, 55414, USA. \EMAIL{wcwj@umn.edu}}
\AUTHOR{Zikai Xiong}
\AFF{MIT Operations Research Center, Cambridge, MA, 02139, USA. \EMAIL{zikai@mit.edu}}
\AUTHOR{Yinyu Ye}
\AFF{Department of Management Science and Engineering, Stanford University, Stanford, CA, 94305, USA. \EMAIL{yyye@stanford.edu}}
} 

\ABSTRACT{%
Identifying optimal basic feasible solutions to linear programming problems is a critical task for mixed integer programming and other applications. The crossover method, which aims at deriving an optimal extreme point from a suboptimal solution (the output of a starting method such as interior-point methods or first-order methods), is crucial in this process. 
This method, compared to the starting method, frequently represents the primary computational bottleneck in practical applications. We propose approaches to overcome this bottleneck by  exploiting problem characteristics and implementing customized strategies. For problems arising from network applications and exhibiting network structures,  we take advantage of the graph structure of the problem and the tree structure of the optimal solutions. Based on these structures, we propose a tree-based crossover method, aiming to recovering basic solutions by identifying nearby spanning tree structures. For general linear programs, we propose recovering an optimal basic solution by identifying the optimal face and employing controlled perturbations based on the suboptimal solution provided by interior-point methods. We prove that an optimal solution for the perturbed problem is an extreme point, and its objective value is at least as good as that of the initial interior point solution. Computational experiments show significant speed-ups achieved by our methods compared to state-of-the-art commercial solvers on classical linear programming problem benchmarks, network flow problem benchmarks, and optimal transport problems.
}


\KEYWORDS{linear programming; crossover; optimal transport; network flow problem; first-order method; interior-point method}

    \maketitle

\section{Introduction} \label{sec: Introduction}
Linear programming is a fundamental problem and important tool in operations research, widely used in various applications like transportation, scheduling, inventory control, and revenue management \citep{bowman1956production,charnes1954stepping,hanssmann1960linear,liu2008choice}. 
A linear programming instance is usually called a linear optimization problem, also known as a linear program or simply ``LP.''
In these problems, \textit{basic feasible solutions} (BFS) -- also known as \textit{corner points} -- correspond identically to \textit{vertices} and \textit{extreme points} of the feasible set. These solutions are particularly important for several reasons: Firstly, to warm start the simplex method, a BFS is required.
Secondly, a BFS, compared to other feasible points, has an inclusive-maximal set of active constraints.
This feature is especially valuable for certain applications, such as discrete optimization problems. Additionally, in integer programming, finding the optimal BFS for LP relaxation problems is essential.
In certain instances, such as network flow problems with integral right-hand side vectors, an optimal BFS is inherently integral. In such cases, finding an optimal BFS for the LP relaxation is sufficient to solve the  integer program.

\subsection{Challenges for Current Crossover Approaches}
Crossover allows us to combine the advantages of both simplex methods and interior-point methods -- the two classic algorithms that have driven advancements in operations research for decades. Most LP solvers use one or both of these methods as the main algorithm. Simplex methods perform pivot iterations through basic solutions and provide an optimal BFS at the end, but they may be outperformed by interior-point methods in solving some problems.
Although interior-point methods converge faster in both theory and practice in solving these problems, they return an interior-point solution rather than a BFS. By employing crossover, it is possible to leverage the speed of interior-point methods for such problems while still obtaining a BFS.

\paragraph{Crossover in LP solvers.}

Crossover has been an essential research topic in conjunction with the studies on interior-point methods. \citet{megiddo1991finding} demonstrated that, given a complementary primal-dual pair of optimal solutions, an optimal BFS could be identified after at most $n$ pivots, where $n$ is the number of variables. Each pivot pushes one variable to the upper or lower bound while maintaining feasibility and optimality of the solution. Although the theory is elegant, in practice, crossover may not satisfy the guarantees of the theoretical results, since it begins with a solution that is only approximately complementary and feasible. To address this issue, \citet{bixby1994recovering} proposed an approach that first identified a near-optimal candidate basis, which might even be infeasible,  by ranking the magnitude of variables' distances to the bounds. They then adapted Wolfe's piecewise-linear phase I algorithm \citep{wolfe1965composite}. Upon reaching a primal feasible solution, they applied \textsc{Cplex}'s phase-II algorithm to solve for the optimal solution. Indeed, Megiddo's crossover is a special case of Bixby and Saltzman's method, which was then included in \textsc{Cplex} 2.2 \citep{bixby1994commentary}. Furthermore, \citet{andersen1996combining} showed that once the interior-point method's solution was sufficiently accurate, it formed an optimal complementary primal-dual pair for a modified LP, whose optimal basis was also an optimal basis for the original problem. In this case, directly applying Megiddo's approach for the modified problem yields an optimal basis for the original problem.  However, \citet{andersen1999exploiting} noted that ensuring the accuracy of the interior-point method remained challenging, and the candidate basis obtained could still be infeasible; thus, in practice, more pivots are needed to arrive at an optimal feasible basis. This approach was subsequently implemented by MOSEK \citep{andersen2000mosek}.

In general, the procedure of solving an LP typically consists of the following three phases, with the second and third phases together referred to as crossover. 

\begin{itemize}
\item \textit{Starting method}: Obtain a suboptimal solution using a ``starting method''. For example, when the main algorithm is set as an interior-point method (also called ``barrier method'' sometimes) in the solver, then the interior-point method serves as the starting method. 
\item \textit{Basis Identification}: Identify a candidate basis from the suboptimal solution and then perform pivots to obtain a nearby primal/dual feasible basis.
\item \textit{Reoptimization}: Adopt primal/dual simplex methods starting from the obtained primal/dual feasible basis until reaching an optimal BFS.
\end{itemize}

The commercial solver \textsc{Gurobi} adopts a similar strategy but it names ``push phase'' the task of performing pivots to obtain a nearby primal/dual feasible basis, while ``cleanup'' corresponds to ``reoptimization'' \citep{Maes2023Initial,gurobi}.
Since the distance of the candidate basis to feasible or optimal bases is unknown, crossover's running time does not have a reliable theoretical estimate. Note that the above three phases differ from the two phases of some simplex methods, where the first phase computes a feasible solution and the second phase computes an optimal solution.

In spite of crossover's significance for practical performance, the literature on this topic is not as extensive as one might expect.  \citet{berkelaar1999basis} extended crossover to quadratic programming. \citet{glavelis2018improving} applied the idea of crossover to accelerate the exterior-point simplex methods. \citet{schork2020implementation} proposed an interior-point method with basis preconditioning and its associated crossover approach. \citet{galabova2020idiot} revealed that the open-source solver Clp \citep{clp} implements a crash method (whose idea comes from \citet{forrest1992steepest}) to increase the sparsity of the given solution via a penalty method. 
\citet{el1994study} and \citet{ye1992finite} studied how to identify the optimal face (or the active constraints at optima); in the special case of non-degeneracy, it is equivalent to identifying the optimal BFS. 

\paragraph{Challenges for modern crossover.}
In many recent real LP applications, crossover has become a crucial component, on average taking over a quarter of the overall running time, as \citet{Maes2023Initial} presented. Given the lack of guarantees on the quality of the candidate basis, the number of pivots can often be enormous, especially for huge-scale problems. Moreover, these large-scale instances are becoming increasingly common in real-world applications.

Crossover also faces challenges from new methods used as the starting method. For general-purpose LPs, many first-order methods have been proposed recently to address huge-scale applications \citep{o2016conic,lin2020admm,deng2022new,applegate2021practical,li2020asymptotically,lu2023cupdlp,wang2023linear}. Although first-order methods can handle larger-scale LPs, obtaining a reliable and highly accurate solution remains a recognized challenge.  This challenge emphasizes the need for further improvements of crossover methods. Furthermore, for LPs with special structures, fast problem-focused first-order algorithms have been proposed \citep{cuturi2013sinkhorn,ge2019interior,altschuler2017near,gao2021boosting}. 
These methods usually do not obtain basic solutions, and the traditional crossover method, designed for general-purpose LPs, cannot utilize the special problem structure. This limitation makes crossover a computational bottleneck in solving such structured problems.

Furthermore, in some cases, obtaining a nearly-optimal BFS within a certain given tolerance is sufficient, which may benefit from early stopping the starting method and dropping the phase of reoptimization. 
However, both the running time of the basis identification phase and the optimality error of the BFS obtained rely on the quality of the solution provided by the starting method. Sometimes, early stopping the starting method cannot guarantee benefits and may even have harmful effects. To deal with this issue, \citet{mehrotra1991finding} proposed a controlled perturbation strategy on the objective vector so that the perturbed problem only has a unique optimal solution. When the optimal solution is unique,  identifying the optimal basis becomes easier, as noted by \citep{tapia1991optimal}. Although Mehrotra's perturbation strategy has demonstrated promising results in practical tests, there is no guarantee that the BFS obtained will have a better objective than the solution provided by the starting method.

In this paper, we propose two novel crossover approaches to address these issues. First, for the \textit{minimum cost flow} (MCF) problem and other problems with a network structure, we present a crossover approach that exploits the graph structure of the problem and the tree structure of the BFS. This approach enables efficient identification of an optimal or near-optimal basis, even when starting from a low-accuracy solution provided by the starting method.
Second, for general LPs, we propose an alternative crossover approach based on detecting the optimal face and applying a minor perturbation. This approach aims to obtain a near-optimal BFS with an objective value at least as good as the solution given by an interior-point starting method.

\subsection{Minimum Cost Flow Problem}
As a prominent and important case of LP, the minimum cost flow (MCF) problem possesses a notable property: its incidence matrix is \textit{totally unimodular}. This means that when the right-hand side vector (in standard form) is integral, every basic solution is also integer-valued. In such cases, if the problem is the LP relaxation of an integer program, one gets an integer solution for free by solving the LP to an optimal BFS. It is worth mentioning that the complexity of the simplex method can be further reduced by the network simplex method, which requires a number of pivots that is at most polynomial in the number of nodes and arcs \citep{orlin1997polynomial}.
Fast algorithms have also been proposed to deal with some special cases, such as \citep{applegate2021practical} for huge-scale problems, and the Sinkhorn method \citep{cuturi2013sinkhorn} for optimal transport problems. These new methods are usually not designed to yield an optimal solution that is basic feasible, so a subsequent crossover is required to obtain an optimal BFS.

The optimal transport problem, as a class of MCF problems, is an important research topic with applications in mathematics \citep{santambrogio2015optimal}, economics \citep{galichon2018optimal}, and machine learning \citep{nguyen2013convergence, ho2019probabilistic}. Machine learning applications in particular require high calculation efficiency of optimal transport problems. By using entropy regularization, the Sinkhorn algorithm \citep{cuturi2013sinkhorn} has significantly reduced the computational complexity and triggered a series of research works. However, in spite of the increased speed, the solution obtained by these methods is not an optimal BFS.

\subsection{Our Contributions}

In this paper, for simplicity of notation, we consider LPs in the following standard form: 
\begin{equation}\label{pro: General LP}
    \min_{x\in\mathbb{R}^n}\,\,\,\,  c^\top x  \quad \operatorname{s.t.}\,\, Ax = b,~  x \ge 0,
\end{equation}
where $A\in\mathbb{R}^{m\times n}$, $c \in \mathbb{R}^n$ and $b\in\mathbb{R}^m$.
\begin{itemize}
    \item 
     For problems with an inherent network structure, we propose two basis identification strategies and their corresponding crossover approaches. Based on the spanning tree characteristics of basic feasible solutions, we develop a tree-based basis identification method. For problems with a large number of arcs, we develop a column-based strategy. These basis identification phases do not require any knowledge of the dual solution but exploit the graph structure of the problem. We conduct experiments with MCF problems on public network flow datasets and test optimal transport problems on the MNIST dataset. Compared with commercial solvers' crossover, our crossover method exhibits significant speed-ups, especially when the given starting method solution is only of low accuracy. Furthermore, we also describe how the combination of first-order methods, such as the Sinkhorn algorithm, and our crossover method, has advantages in practice.
    
    \item For general LPs, we propose a crossover approach that first identifies an approximate optimal face using the primal-dual solution pair and then finds a high-quality near-optimal BFS by applying a minor perturbation. Depending on the requirements of applications, the reoptimization phase can be carried out via the primal or dual simplex method.     
    Experiments on LP benchmark problems reveal that our crossover method can notably enhance the performance of state-of-the-art commercial solvers, particularly on currently challenging problems.

\end{itemize}

\subsection{Outline}
In \Cref{sec: preliminary}, we introduce some preliminaries and background. 
In \Cref{sec: Network Crossover Method}, we propose a network structure based crossover method for MCF problems and its extension to general LPs. 
In \Cref{sec: Perturbation Crossover Method}, we propose a perturbation crossover method to accelerate the crossover procedure in general-purpose LP solvers. The computational results are presented in \Cref{sec: Numerical Experiments}.

\subsection{Notation}
We use $[n]$ to denote the set $\{1,2,\dots,n\}$.  For a set $\mathcal{X}$, the notation $\vert \mathcal{X}\vert$ denotes the cardinality of $\mathcal{X}$. For $\mathcal{X}\subseteq [n]$ and vector $u\in\mathbb{R}^n$, $u_{\mathcal{X}}$ is the $\vert \mathcal{X}\vert$-vector constituted by the components of $u$ with indices in $\mathcal{X}$. 
For a matrix $X$ and scalar $x$, $X \ge x$ means each component of $X$ is greater than or equal to $x$. For a matrix $A = (A)_{ij}$, the positive part and negative part of $A$ are $A^+$ and $A^-$, $(A^+)_{ij} := \max\{(A)_{ij}, ~0\},$ and $(A^-)_{ij} := \max\{-(A)_{ij},0\}.$ $A_{i\cdot}$, $A_{\cdot j}$ and $A_{ij}$ denote the $i$-th row, $j$-th column and the component in the $i$-th row and the $j$-th column of $A$, respectively. 
For a graph $(\mathcal{N}, \mathcal{A})$ and $i\in \mathcal{N}$, $\mathcal{O}(i),\mathcal{I}(i)$ denote the node sets connecting to i, specifically, $\mathcal{O}(i):=\{j: (i,j)\in \mathcal{A}\}, \mathcal{I}(i):=\{j: (j,i)\in \mathcal{A}\}$. For any two sets $\mathcal{X}$ and $\mathcal{Y}$, we use $\mathcal{X} \setminus \mathcal{Y}$ to denote $\{x: x\in \mathcal{X} \text{ but }x\notin \mathcal{Y}\}$. 
We denote the Euclidean norm using $\|\cdot\|$.
For any set $\mathcal{X}$, we use $P_{\mathcal{X}}(\cdot)$ to denote the projection onto $\mathcal{X}$, i.e., $P_{\mathcal{X}}(z):=\arg\min_{\hat{z}\in\mathcal{X}}\|z - \hat{z}\|$. For a matrix $A \in \mathbb{R}^{n\times n}$, we use $A^\dag$ to denote the Moore–Penrose inverse of $A$. Moreover, we use $\operatorname{Null}(A)$ and $\operatorname{Im}(A)$ to denote the null space and image space of the matrix $A$, respectively.

\section{Preliminaries}\label{sec: preliminary}

\paragraph{MCF problem.}

Let $\mathcal{G} = (\mathcal{N},\mathcal{A})$ be a directed graph consisting of nodes $\mathcal{N}$ and arcs $\mathcal{A}$. Each node $i\in \mathcal{N}$ is associated with a signed supply value $b_i$. Each arc $(i,j)\in \mathcal{A}$ is associated with a capacity $u_{ij}$ (possibly infinite) and a cost $c_{ij}$ per unit of flow. 
Then the MCF problem is defined as
\begin{equation}\label{pro: General MCF Problem}
\begin{array}{ll}
\min_f & \sum_{(i,j)\in\mathcal{A}}c_{ij}f_{ij} \\
\operatorname{s.t.}& b_i + \sum_{j\in \mathcal{I}(i)}f_{ji} = \sum_{j\in \mathcal{O}(i)}f_{ij},~ \forall~ i \in \mathcal{N}, \text{ and }0\le f_{ij} \le u_{ij},\ \forall  \ (i,j)\in \mathcal{A},
\end{array}
\end{equation}
where $\mathcal{O}(i) := \{j: (i,j)\in \mathcal{A}\}$, $\mathcal{I}(i) := \{j: (j,i)\in \mathcal{A}\}$.

\paragraph{Optimal transport problem.} The optimal transport problem is a special equivalent form of the general MCF problems. Following the definition of the general MCF problem, the nodes $\mathcal{N}$ can be divided into two groups, suppliers $\mathcal{S}$ and consumers $\mathcal{C}$. And arcs $\mathcal{A}$ pairwisely connect suppliers and consumers, in the direction from the former to the latter, i.e., $\mathcal{A} = \mathcal{S}\times \mathcal{C}:= \{(i,j):i\in \mathcal{S},j\in\mathcal{C}\}$. For any supplier $i$ and consumer $j$, $s_i$ and $d_j$ respectively denote the capacity of the supply and demand.
Then an optimal transport problem can be rewritten as follows:
\begin{equation}\label{pro: General OT Problem}
\begin{array}{lll} \operatornamewithlimits{min}_{f} & \sum_{(i,j)\in\mathcal{S}\times \mathcal{C}} c_{ij}f_{ij} \\
\operatorname{s.t.}& \sum_{j\in \mathcal{C}}f_{ij} = s_i, ~ \forall ~ i \in \mathcal{S}, \ \sum_{i\in \mathcal{S}}f_{ij} = d_j, ~ \forall j \in \mathcal{C}, \text{ and } f_{ij} \ge 0, ~ \forall ~ (i,j)\in \mathcal{S}\times \mathcal{C}
\end{array}
\end{equation}

\paragraph{Optimal Face.}
For the general primal problem \eqref{pro: General LP}, the  dual problem is
\begin{equation}\label{eq: dual LP}
\begin{array}{l}
   \operatornamewithlimits{max}_{y\in \mathbb{R}^m,s\in\mathbb{R}^n}  ~ b^\top y \quad
   \operatorname{s.t.} ~ A^\top y + s = c,~~ s\ge 0
\end{array}
\end{equation}
Specifically, any feasible $(y,s)$ satisfies $s = c - A^\top y$ so we say a slack $s\ge 0$ is feasible whenever there exists $y$ such that $A^\top y + s = c$. It can be further proven that for any feasible $s$, all the corresponding $y$ must share the same objective value, which is a linear function of $s$. This is because for any $y$ satisfying $A^\top y + s =c $, it holds that $b^\top y = (A\tilde{x})^\top y = \tilde{x}^\top (c - s)$ for any $\tilde{x}$ such that $A\tilde{x} = b$.
If we replace $y$ with its corresponding $s$, then the problem in $s$ becomes symmetric with the primal problem, as the feasible set becomes the intersection of the affine subspace $\{s: \exists \ y\in \mathbb{R}^m \text{ such that }A^\top y + s = c\}$ and $\mathbb{R}^n_+$, and the objective becomes a linear function of $s$.  This idea of analyzing the dual in a  way symmetric to the primal, to the best of our knowledge, was first proposed by \cite{todd1990centered}, and also used in classic textbooks such as \citep{renegar2001mathematical}.

For \eqref{pro: General LP} and \eqref{eq: dual LP}, there exists at least one optimal solution pair $(x^\star,s^\star)$ that is strictly complementary \citep{goldman1956theory}, i.e., $\mathcal{P}(x^\star)\cap \mathcal{P}(s^\star) = \emptyset$ and $\mathcal{P}(x^\star)\cup \mathcal{P}(s^\star) = [n]$,
in which $ \mathcal{P}(z):= \{i:z_i>0\}$ denotes the set of indices for the strictly positive components of $z$.
Moreover, if we denote $\mathcal{P}(x^\star)$ as $\mathcal{P}^\star$ , and denote $\mathcal{P}(s^\star)$ as $\mathcal{D}^\star$, then $\{\mathcal{P}^\star,~\mathcal{D}^\star\}$ is called the \textit{optimal partition}. For any optimal solutions $x$ and $s$, $\mathcal{P}(x) \subseteq \mathcal{P}^\star,$ and $ \mathcal{P}(s) \subseteq \mathcal{D}^\star$.
With the above definitions, the optimal face of the primal problem \eqref{pro: General LP} and dual problem \eqref{eq: dual LP} can be rewritten as
\begin{equation}\label{eq def theta}
\begin{array}{ll}
&\Theta_p:=\{x: Ax = b, ~x_{\mathcal{P}^\star} \ge 0 , \text{and } x_{\mathcal{D}^\star} = 0\}, \\
& \Theta_d :=\{(y,s):A^\top y + s = c,~s_{\mathcal{D}^\star}\ge0, \text{ and }s_{\mathcal{P}^\star}  = 0\}.
\end{array}
\end{equation}


\paragraph{Column generation method.}

The column generation method is an approach for large-scale LPs. Instead of directly solving the original problem \eqref{pro: General LP}, the column generation method involves a sequence of master iterations. In each iteration, a collection of columns $A_{\cdot i}$, $i\in \Lambda$, is formed and the following restricted problem is solved:
\begin{equation}\label{pro: column generation restricted problem}
    LP(\Lambda):\quad \min_{x\in\mathbb{R}^n}~  \sum_{i\in \Lambda}c_i x_i ~~ \operatorname{s.t.}~ \sum_{i\in \Lambda}A_{\cdot i}x_i = b,~ x_{[n]\setminus \Lambda}=0,~ x_{\Lambda} \ge 0.
\end{equation}
Refer to \Cref{alg: column generation method framework} for the general framework. See \cite{bertsimas1997introduction} for specific rules on updating $\Lambda_k$.  The collection $\Lambda_{k+1}$ contains the basis of the solution $x^k$ from $LP(\Lambda_k)$, so a simplex method for $LP(\Lambda_{k+1})$ can easily warm start from $x^k$.

\begin{algorithm}[H]
\caption{General Column Generation Method}
\label{alg: column generation method framework}
\begin{algorithmic}[1]
\STATE Iteration counter $k=1$, and initial set of columns $\Lambda_1$.
\WHILE{the reduced cost $\bar{c}\ngeq 0$}
    \STATE \textit{Update column set:} Update $\Lambda_{k}$ according to some rules.
    \STATE \textit{Solve the restricted problem:} Solve $LP(\Lambda_k)$ and obtain an optimal basis $B_k$.
    \STATE \textit{Compute reduced cost:} $\bar{c}\gets c-c^\top_{B_k}B_k^{-1}A$, and $k \gets k+1$.
\ENDWHILE
\end{algorithmic}
\end{algorithm}

\section{Network Crossover Method}\label{sec: Network Crossover Method}
As mentioned in \Cref{sec: Introduction}, traditional crossover methods select the candidate basis by directly ranking the magnitude of variables or comparing ratios between variables. However, these methods do not take the graph structure of the problem into account.

\begin{figure}[htbp]
    \centering
    \includegraphics[width=0.92\linewidth]{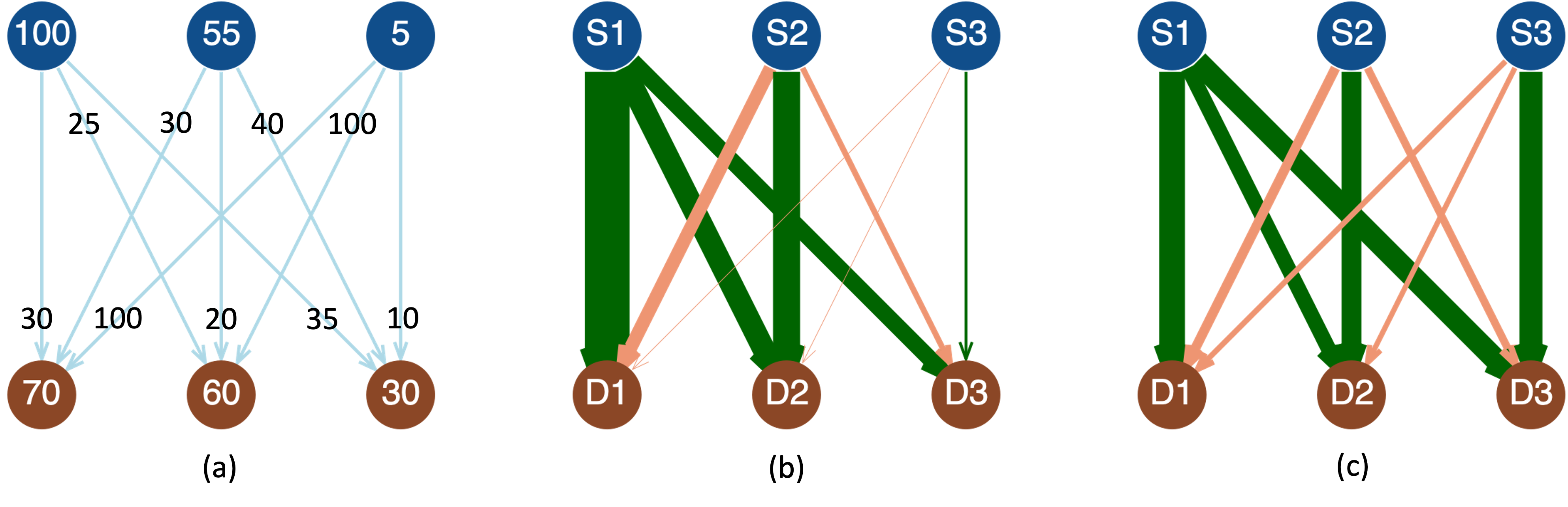}
    \caption{An example of an optimal transport problem, with supply nodes in blue and demand nodes in brown. Figure (a) presents the problem coefficients of supplies, demands, and transportation costs per unit. Figure (b) shows the relative values of an interior solution as the line width of arcs. Figure (c) shows the relative flow ratios defined in \eqref{def: MCF flow ratio} as the line width of arcs.}
    \label{fig:flowratioexample}
\end{figure}

We first use an optimal transport example to show that ranking the magnitude of variables may introduce errors.
\Cref{fig:flowratioexample} illustrates an example of a small optimal transport problem with supply nodes in blue and demand nodes in brown. The supplies, demands, and transportation costs per unit are displayed in  \Cref{fig:flowratioexample}a.
Suppose that we are given a central-path solution generated by a starting method, such as a typical path-following interior-point method, and this solution has 
an objective value that is 7\% higher than the optimal objective value. 
Then the relative variable values of this solution are represented in the arc widths of \Cref{fig:flowratioexample}b.
The arcs in dark green denote the variables in the optimal basis and the arcs in light red denote the variables not in the optimal basis. 
From the problem setting, the arc $(S3,D3)$ should be included in the optimal basis; however, its width is not as great as that of arcs $(S2,D1)$ or $(S2,D3)$. Indeed, by taking into account the graph structure, since every basis corresponds to a spanning tree, the optimal basis must contain an arc that originates from node $S3$. 
Therefore, in this section, we will define the \textit{flow ratio} (whose value is represented by the arc widths in \Cref{fig:flowratioexample}c) to depict the graph structure and use it for crossover.  

Before we define the flow ratio, we first transform the MCF problem into an uncapacitated problem based on the value of the variables at a given suboptimal solution.
Assume that we have obtained a feasible suboptimal solution $\tilde{f}$ for the problem \eqref{pro: General MCF Problem}, denoted as $MCF(b,c,u)$, using a starting method. In this case, there exists a transformation of the problem variables such that each $\tilde{f}_{ij}$ is closer to the lower bound $0$ than the capacity upper bound $u_{ij}$. Specifically, let $\mathcal{L} := \{(i,j)\in\mathcal{A}: \tilde{f}_{ij} - 0 < u_{ij} - \tilde{f}_{ij}\}$ and $\mathcal{U}:= \mathcal{A} \setminus \mathcal{L}$. Our requirement is then satisfied by replacing variables $f_{\mathcal{U}}$ in $MCF(b,c,u)$ with $u_{\mathcal{U}} - f_{\mathcal{U}}$. We denote the problem after the variable transformation as $MCF(\bar{b},\bar{c},\bar{u})$.
Note that if $\tilde{f}$ is already an optimal solution of $MCF(b,c,u)$, then the corresponding $MCF(\bar{b},\bar{c},\bar{u})$ has an optimal solution that is never active in any constraint of the capacity upper bound. Moreover, when $\tilde{f}$ is close enough to the optimal solutions (for example, when the $l_\infty$ norm distance is smaller than $\frac{1}{2}\min_{i}|\bar{u}_i|$), $MCF(\bar{b},\bar{c},\bar{u})$ shares the same optimal solutions with the uncapacitated MCF problem $MCF(\bar{b},\bar{c},\infty)$. Also, any optimal BFS of $MCF(\bar{b},\bar{c},\infty)$, after being transformed back to the original MCF problem, remains an optimal BFS of the original problem. Consequently, to identify a near-optimal basis for $MCF(b,c,u)$, we can instead identify a near-optimal basis for $MCF(\bar{b},\bar{c},\infty)$ by constructing a nearby basis for the transformation of the given flow vector $\tilde{f}$.

In uncapacitated MCF problems, a solution $f$ is a basic (feasible) solution if and only if it is a (feasible) \textit{tree solution} \citep{ahuja1993network}. The tree solutions in uncapacitated MCF problems are those whose nonzero components can form a spanning tree without considering the direction of flows. For a nondegenerate BFS of the uncapacitated MCF problem in a connected graph, there are exactly $\vert \mathcal{N}\vert -1$ strictly positive flows (positive components of the BFS), connecting each node in the graph and corresponding to a spanning tree in $\mathcal{G}$.
Taking into account the tree structure of basic solutions, we will select $|\mathcal{N}| -1$ arcs in a specific order. For a nonbasic suboptimal feasible solution, the arc that delivers the most flow to a single node has a higher likelihood and priority of being in the basis compared to other flows.
Therefore, we use the following \textit{flow ratio} to measure the potential of an arc to be part of an optimal basis.

Given any flow $f$ of a directed graph $\mathcal{G} = (\mathcal{N},\mathcal{A})$, for each arc $(i,j)$ in $\mathcal{A}$, the \textit{flow ratio} is defined as the maximum proportion of the flow that $f_{ij}$ serves over the total flow that passes node $i$ or $j$, i.e.,
\begin{equation}\label{def: MCF flow ratio}
    r_{ij}:= \max\left\{
    \frac{f_{ij}}{\sum_{l\in \mathcal{O}(i)}f_{il}+ \sum_{l\in \mathcal{I}(i)} f_{li}}, \
    \frac{f_{ij}}{\sum_{l\in \mathcal{O}(j)}f_{jl}+ \sum_{l\in \mathcal{I}(j)} f_{lj}}
    \right\}\ .
\end{equation}
The flow ratio measures the proportion of the flow in a specific arc compared to the total flow in all the conjoint arcs. All optimal solutions are convex combinations of optimal basic feasible solutions, so for any method generating iterates that converge to optimal solutions, the flow ratio of the arcs that are not part of any optimal basis would converge to zero. Unlike traditional indicators that only consider the magnitude of the arc, the flow ratio takes into account the connected nodes and all adjacent arcs.  For instance, in \Cref{fig:flowratioexample}c, the linewidths represent the relative flow ratios, which align better with the optimal basis than when solely using variable values in \Cref{fig:flowratioexample}b. 
Same with the traditional crossover methods, there is no strong theoretical guarantee for these indicators,  but we will provide two powerful crossover methods based on them.  

Recall that, after the starting method, crossover is composed of two phases, namely basis identification and reoptimization. 
Let the arcs be ranked in descending order of their corresponding flow ratios and form the array $(s_1,s_2,\dots, s_{|\mathcal{A}|})$ of $\mathcal{A}$. Based on this array, we propose two crossover methods for the basis identification phase.

\subsection{Tree-based Basis Identification (\textsc{Tree-BI})}


Similarly to the standard crossover, the basis identification first identifies a candidate basis and then obtains a nearby feasible basis as a starting point for reoptimization. 
Since the BFS of any MCF problem has to be a feasible tree solution, and a tree solution corresponds to a spanning tree in the graph, candidate bases can be identified from the construction of a spanning tree.
If the given suboptimal solution is sufficiently accurate and the optimal solution is unique, such a spanning tree corresponds exactly to the optimal BFS. This spanning tree problem can be efficiently solved using combinatorial algorithms for minimum spanning tree problems \citep{prim1957shortest,dijkstra1959note,chazelle2000minimum}.
When the suboptimal solution is not sufficiently accurate, the basis obtained may still be infeasible. In such cases, in a fashion similar to traditional crossover methods, a piecewise-linear phase I simplex method is then required to find a nearby feasible basis.

For the special case of optimal transport problems, the following proposition enables us to design an efficient phase I method for obtaining a feasible basis:
\begin{proposition}\label{pro: pushphase}
For the optimal transport problem \eqref{pro: General OT Problem}, suppose that the basic solution/flow $f$ is infeasible, say $f_{ij}<0$, for $(i,j)\in\mathcal{A}$. Then, there exist $j', i'\in\mathcal{N}$ such that $(i,j')$, $(i',j)$ are also arcs in $\mathcal{A}$ and $f_{ij'}>0$, $f_{i'j}>0$. Moreover, for any such $i', j'$, $(i',j')$ is also an arc in $\mathcal{A}$ and $f_{i'j'}=0$.
\end{proposition}
\proof{Proof.}
    Since $(i,j)\in \mathcal{A}$, we have $i\in \mathcal{S}$ and $j\in\mathcal{C}$. Moreover, $\sum_{l\in\mathcal{C}}f_{il} = s_i \ge 0$ and $\sum_{l\in S}f_{lj} = d_j \ge 0$ due to the definition of optimal transport problems. Noting that $f_{ij}< 0$, there must exist $f_{ij'}>0$, $j'\in \mathcal{C}$ and $f_{i'j}>0$, $i'\in \mathcal{S}$. Furthermore, because $i'\in\mathcal{S}$, $j'\in \mathcal{C}$ and $\mathcal{A} = \mathcal{S} \times \mathcal{C}$, then $(i',j')$ is an arc in $\mathcal{A}$. 
    Since the subgraph corresponding to the nonzeros in $f$ is a tree, it must contain no cycle. In particular, there must be zero flow in $f_{ij}, f_{ij'}, f_{i'j},f_{i'j'}$, and thus $f_{i'j'} = 0$.
\endproof
Using this result, we can eliminate the negative flows of the candidate basic solution, denoted as $f_{ij}$, by repeatedly executing the following steps:
\begin{enumerate}
    \item {Find adjacent positive flows:} $j' \leftarrow \textup{arg max}_{l\in\mathcal{C}} f_{il}$ \text{, } $i' \leftarrow \textup{arg max}_{k\in\mathcal{S}} f_{kj}$;
    \item {Adjust the flows:} Let $\theta \leftarrow \min\{-f_{ij}, f_{i'j}, f_{ij'}\}$ and then update the flow: $f_{ij} \leftarrow f_{ij} + \theta$,
		$f_{ij'} \leftarrow f_{ij'} - \theta$,
		$f_{i'j} \leftarrow f_{i'j} - \theta$,
		$f_{i'j'} \leftarrow \theta$.
\end{enumerate}
Note that this procedure does not influence any other negative flows but decreases the negativity of $f_{ij}$, so repeating it must lead to a feasible basic flow in the end. Furthermore, the $j'$ and $i'$ found in step 1 do not have to be the maximum in practice.

\subsection{Column-based Basis Identification (\textsc{Col-BI})}
Although \textsc{Tree-BI} can identify a candidate basis by solving a minimum-spanning tree problem using fast combinatorial algorithms, there is no guarantee that the obtained nearby feasible BFS will remain close to the variables with high flow ratios after the additional pivots. To address this issue, we propose an alternative method that utilizes column generation to identify a  candidate basis, where high flow ratios are maintained at the basic variables. We refer to this method as \textsc{Col-BI}.

Let $E \in \mathbb{R}^{|\mathcal{N}|\times|\mathcal{A}|}$ be the node-arc incidence matrix of the network $\mathcal{G}$, defined as follows: $E_{ij} := 1$ if the $j$-th arc leaves node $i$, $E_{ij} := -1$ if the $j$-th arc enters node $i$, and $E_{ij} := 0$ for all other cases.
Then, the MCF problem \eqref{pro: General MCF Problem} is equivalent to the \textit{big-M method} formulation on the graph with one additional root node and artificial flows in arcs between each node in $\mathcal{N}$ and the root node. This formulation is a common approach for initializing the network simplex method. If we use indices $1,\dots, |\mathcal{A}|$ to denote the arcs and indices $|\mathcal{A}| + 1 ,\dots,|\mathcal{A}| + |\mathcal{N}|$ to denote the artificial arcs, the corresponding LP formulation is as follows:
\begin{equation}\label{pro: Modified LP with artifical variables}
    \begin{array}{lll}
    \min_{\bar{f}:=(f;\hat{f}) \in \mathbb{R}^{|\mathcal{N}|+|\mathcal{A}|}} & c^\top f   +  M e^\top \hat{f} &\\
    \operatorname{s.t.}& E f + D \hat{f} = b, \ e^\top D \hat{f} = e^\top b \\
    & 0\le f \le u,\, \hat{f}\ge0
    \end{array}
\end{equation}
where $e$ is the all-one vector, $D$ is a diagonal matrix in $\mathbb{R}^{|\mathcal{N}|\times |\mathcal{N}|}$ with the diagonal entries being $\{\operatorname{sgn}(b_i)\}_{i=1}^{|\mathcal{N}|}$, and $M$ is a large positive scalar, such as $|\mathcal{A}| \cdot \|c\|_{\infty}$. Problem \eqref{pro: Modified LP with artifical variables} is equivalent to \eqref{pro: General MCF Problem} but has an obvious BFS in which $f$ is a zero vector and $\hat{f}$ is a nonnegative vector with components equal to the absolute values of $b$'s entries. We introduce problem \eqref{pro: Modified LP with artifical variables} because it can be used for the column generation method in Algorithm \ref{alg: column generation method framework}. The column sets are then generated as $\Lambda_0:=\{|\mathcal{A}|+1,\dots,|\mathcal{A}|+|\mathcal{N}|\}$ and 
\begin{equation} \label{eq: generate series of subproblem}
    \Lambda_{k+1}  \gets  \bar{\Lambda}_k \cup \big( \Lambda_0  \setminus \{i:\bar{f}^{k}_{i}\text{ is nonbasic and artificial}\} \big) \ , 
\end{equation}
where $\bar{f}^{k}$ is the solution of the $k$-th problem $LP(\Lambda_k)$,  and $\bar{\Lambda}_k$ is the column sets that could be priced at iteration $k$. We choose $\bar{\Lambda}_k$ according to the arc order $(s_1,\dots,s_{|\mathcal{A}|})$ to represent the columns that may contain the optimal basis. We set the number of components in $\bar{\Lambda}_k$ to be $\theta_k$ for a certain monotonically increasing series $(\theta_k)_{k\in\mathbb{N}}$, and define $\bar{\Lambda}_k$ as $\{i_{s_1},i_{s_2},\dots,i_{s_{\theta_k}}\}$. Here the $i_{s_j}$ is the index of the column in the formulation \eqref{pro: Modified LP with artifical variables} corresponding to the arc $s_j$ in the graph. 
In this way, the column-generation method starts from the given artificial basis and progressively incorporates additional columns with high flow ratios. Concurrently, it eliminates the artificial arcs of zero variables. The restricted problem is a smaller MCF problem, and it can be solved by the network simplex method more efficiently.  Once all artificial variables reach zero, the solution of the restricted subproblem emerges as a BFS for the original problem.

\subsection{Reoptimization Phase}

After obtaining a feasible BFS from either \textsc{Tree-BI} or \textsc{Col-BI}, we can proceed to the reoptimization phase to obtain an optimal BFS. In this phase, the column generation algorithm \ref{alg: column generation method framework} on the original MCF problem is still used to give higher priorities to the columns with higher flow ratios when doing pivots. The sequence $(\Lambda_k)_{k\in\mathbb{N}}$ is initialized by $\Lambda_0$ equal to the starting basis and then generated by involving the columns with negative reduced cost or high flow ratios, i.e., $
    \Lambda_{k+1} \gets \Lambda
    _k \cup \{i:\bar{c}_i < 0\} \cup\{i_{s_1},i_{s_2},\dots,i_{s_{\theta_k}}\},$
until $\bar{c} \ge 0$, where $(\theta_k)_{k\in\mathbb{N}}$ is any monotonously increasing integer series.

\subsection{Extension to General LP Problems}\label{sec: Extension to General LP Problems}


The MCF problem is only a special case of LPs, and the definition of the flow ratio can be extended to general LPs. For the general primal LP \eqref{pro: General LP}, we define the general flow ratio as, for each $i\in [n]$, 
\begin{equation}\label{def: general LP flow ratio}
    r_i := \max\left\{\frac{|A_{ki}x_i|}{\sum_{j = 1}^n |A_{kj}x_j| }: k=1,\dots m, \text{ and }A_{ki}x_i \neq 0\right\} 
\end{equation}
to represent the maximal proportion that $x_i$ serves in a certain row $k$'s sum of $\sum_j |A_{kj}x_j|$. The definition of flow ratio for MCF problems in \eqref{def: MCF flow ratio} is only a special case of \eqref{def: general LP flow ratio}. The  approach of using the general flow ratio \eqref{def: general LP flow ratio} to do crossover for general LPs is exactly the same as the approach for MCF problems, except that the network simplex method is replaced with the standard simplex method. For LPs that are not MCF problems but share similar structures, the flow ratio could still capture both the problem structure and the information of the suboptimal solution.

\section{Perturbation Crossover Method}\label{sec: Perturbation Crossover Method}

The network crossover method in \Cref{sec: Network Crossover Method} is for problems with network structure. In this section, we propose a new perturbation based crossover procedure for cases where the starting method is an interior-point method that returns a highly accurate interior-point solution. This is the most typical scenario for the traditional crossover associated with interior-point methods. This method combines two techniques, optimal face detection~(\Cref{subsec: optimal face detection}) and smart perturbation~(\Cref{subsec: controlled perturbation}), to obtain a high-quality BFS. We prove that this BFS is at least as good as the given interior point solution when the given solution is on the central path.
Before discussing the techniques, we first introduce the following lemma about using perturbations to find a BFS via ensuring only one unique solution exists.
\begin{lemma}\label{lm: perturbation find BFS}
    The set of all $p$ such that the perturbed primal problem
    \begin{equation}\label{pro: general perturbed primal problem}
        \min_x\  (c+p)^\top x \ \ \ \  \operatorname{s.t. } \ Ax = b , \ x \ge 0
    \end{equation}
    has multiple optimal solutions is of zero measure. 
\end{lemma}

\Cref{lm: perturbation find BFS} shows that a random perturbation almost surely ensures the perturbed problem has a unique optimal solution. 
Since the LP with a unique optimal solution must have that optimal solution as a BFS, this result implies that the optimal solution of a randomly perturbed problem is very likely to be a BFS of the original problem. 
Similar perturbations have also been used to get rid of degeneracy in LP instances \citep{charnes1952optimality,megiddo1989varepsilon}. A proof is attached in Appendix \ref{app:proof} for completeness.
Symmetric results also hold for the dual problem by perturbing the right-hand-side vector $b$ to $\bar{b}:=b+p$. 
However, not all random $p$ ensure that $Ax = \bar{b}$ has a feasible solution, especially when $A$ is not full-rank. Therefore, we consider perturbations of the form  $\bar{b}:=b+Ap$ instead. Using a similar proof (omitted because it is analogous), we have the following lemma:

\begin{lemma}\label{lm: perturbation find BFS dual}
    The set of all $p$ such that the perturbed dual problem
    \begin{equation}\label{pro: general perturbed dual problem}
        \max_{y,s} \ (b+ A p)^\top y \ \ \ \ \operatorname{s.t. } \ A^\top y + s = c , \ s \ge 0
    \end{equation}
    has multiple optimal solutions is of zero measure. 
\end{lemma}

\subsection{Optimal Face Detection}\label{subsec: optimal face detection}

\Cref{lm: perturbation find BFS,lm: perturbation find BFS dual} imply that an optimal primal (or dual) BFS can be computed by directly solving an LP whose objective vector is a random vector and whose feasible set is the optimal face $\Theta_p$ or $\Theta_d$, as defined in \eqref{eq def theta}. The good news is that it has been proven that the iterates of interior-point methods provide indicators of the optimal partition, which is used to define the optimal face.
Let $(x^k,s^k)$ be the iterate of interior-point methods. Then, for a class of classic interior-point algorithms, \citet{mehrotra1993finding} proved that when $k$ is sufficiently large, the optimal partition $\{\mathcal{P}^\star,\mathcal{D}^\star\}$ is identical to 
$\{{\mathcal{P}}^k := \{j: x_j^k \ge s_j^k\},{\mathcal{D}}^k := \{j: x_j^k < s_j^k\}\}$. 

However, it is difficult to determine 
the $k$ that is sufficiently large, and the interior-point methods used in solvers may no longer have these theoretical guarantees, so in practice, we may use the following relaxed candidate optimal faces to replace the real optimal faces in the perturbation problem. 
\begin{equation}\label{eq: relaxed identified optimal face}
\begin{array}{ll}
    & \Theta_{p,\gamma}^k:= \{x: Ax = b, ~x\ge 0 , \text{and } x_j = 0 \text{ for } j \notin \mathcal{P}^k_\gamma\} \text{, where }\mathcal{P}^k_\gamma:= \{j: x^k_j \ge \gamma s^k_j\}\\
    & \Theta_{d,\gamma}^k:= \{(y,s): A^\top y + s = c, ~s\ge 0 , \text{and } s_j = 0  \text{ for } j\notin \mathcal{D}^k_\gamma\}\text{, where }\mathcal{D}^k_\gamma := \{j: s^k_j >  \gamma x^k_j\} \ .
\end{array}
\end{equation}
Here, $\gamma \in (0,1]$ is a given relaxation parameter that determines how conservative the candidate optimal faces are. 
If the perturbed problem on $\Theta_{p,\gamma}^k$ or $\Theta_{d,\gamma}^k$ is identified infeasible, the relaxation parameter $\gamma$ is decreased to restrict fewer variables to zero. However, the relaxed candidate optimal face may also be overly conservative and contain nonoptimal solutions, which means directly solving the problem with a random objective vector on the candidate optimal face is not a good approach. Therefore, in the next subsection we introduce how to solve a slightly perturbed problem so that its solution remains nearly optimal for the original problem.

\subsection{Smart Perturbation Based on Interior-Point Method Solutions}\label{subsec: controlled perturbation}
In this subsection, we introduce a smart perturbation technique that leverages the central-path solutions to guarantee the quality of the solutions of the perturbed problems.
For simplicity of analysis, following a classic textbook like \citep{nesterov1994interior}, we assume that both the primal problem \eqref{pro: General LP} and its dual problem are strictly feasible. If the problem is not strictly feasible, practical interior-point methods solve an equivalent reformulation, such as the homogeneous self-dual reformulation, to ensure strict feasibility  \citep{andersen2009homogeneous}.

Practical interior-point methods generate a sequence of solutions $(x^k,y^k,s^k)$ following the central path \citep{nesterov1994interior}. The primal-dual solution pair $(x,y,s)$ is on the central path if it is feasible and all products $x_is_i$ are equal to $\mu$ for a certain $\mu > 0$. We say that $x$ is on the central path if there exists $(y,s)$ such as $(x,y,s)$ is on the central path. We will use $X_k$ and $S_k$ to denote the diagonal matrices with diagonal entries being the components of $x^k$ and $s^k$, respectively. The following theorem will show that if $x^k$ is on the central path, and the size of the perturbation is nicely controlled, then any optimal solutions of the perturbed problem are at least as good as $x^k$ in terms of the corresponding objective value on the original problem. 

\begin{theorem}\label{thm: controlled perturbation for subspace LP}
     Let $x^k$ be on the central path of the LP problem
     \begin{equation}\label{pro: subspace LP}
         \min \ c^\top x, \quad \operatorname{s.t.} \ x \in L':=L + x_0, \ x\in \mathbb{R}^n_+ \ ,
     \end{equation}
     where $L'$ is an affine subspace in $\mathbb{R}^n$, associated with a linear subspace $L$. Let $p$ be any perturbation in $\mathbb{R}^n$ such that \begin{equation}\label{eq:sizeofperturbation}
     \|X_k p\|_2 < \frac{1}{n+2\sqrt{n}+1}\cdot \| P_{X_k^{-1} \cdot L} (X_k c)\|_2 .
     \end{equation}
     Then for the optimal solution $\hat{x}$ of the perturbed problem   \begin{equation}\label{pro: perturbed subspace LP}
         \min \ (c+p)^\top x, \quad \operatorname{s.t.} \ x \in L':=L + x_0, \ x\in \mathbb{R}^n_+ \ ,
     \end{equation}
    its objective value on \cref{pro: subspace LP} is no higher than that of $x^k$, namely, $c^\top \hat{x} \le c^\top x^k$.
\end{theorem}

An instance of \eqref{pro: subspace LP} is the primal problem \Cref{pro: General LP}, where $L = \operatorname{Null}(A)$, and $x_0$ is any solution of the linear system $Ax = b$. For this type of problems, $P_{X_k^{-1} \cdot L} (X_k c)$ has closed form
\begin{equation} \label{eq: projector def}
    P_{X_k^{-1} \cdot L} (X_k c) = (I - X_k A^\top (AX_k^2 A^\top)^\dag A X_k) X_k c.
\end{equation}
\Cref{thm: controlled perturbation for subspace LP} presents  a range of perturbation on the objective vector $c$ such that the solution $\hat{x}$ is no worse than $x^k$ in the objective value. 
Another example of \eqref{pro: subspace LP} is the dual problem \Cref{eq: dual LP} of only the slack $s$, which is in the symmetric formulation of the primal problem, with $L$ being $\operatorname{Im}(A^\top)$ and $x_0$ being $c$.   For the dual problem, \Cref{thm: controlled perturbation for subspace LP} presents a range of the  perturbation on the right-hand-side vector $b$ such that the solution has no worse objective value than interior-point dual solution $(y^k,s^k)$.

Together with \Cref{lm: perturbation find BFS,lm: perturbation find BFS dual}, \Cref{thm: controlled perturbation for subspace LP} implies that, when the given interior-point method solution $x^k$ (or $(y^k,s^k)$) is on the central path and no non-optimal BFS has better objective values of the primal (or dual) problem, then using the random perturbations that satisfy \Cref{eq:sizeofperturbation} on the objective vector $c$ (or the right-hand-side vector $b$) would directly almost surely ensure that the solution of the perturbed problem is an optimal primal (or dual) BFS. 
Even if the interior point solution is not accurate enough, the perturbations that satisfy \Cref{eq:sizeofperturbation} also ensure that the solutions of the perturbed problems will not get worse. We will call the perturbation $p$ that satisfies \Cref{eq:sizeofperturbation} a ``smart perturbation.''

Before proving \Cref{thm: controlled perturbation for subspace LP}, we first recall a classic result about the central-path solutions.
Let the sublevel set of $x^k$ for \eqref{pro: subspace LP} be denoted as $\mathcal{F}_k := \{x: x\in L'\cap \mathbb{R}^n_+,c^\top x \le c^\top x^k\}$, the set of feasible solutions with the objective value no larger than $c^\top x^k$, and then we have:
\begin{lemma}\label{lm: analytic center rounding}
Let $x^k$ be any point on the central path of the problem \eqref{pro: subspace LP} and let $\mathcal{F}_k$ be the corresponding sublevel set. Then for any $\tilde{x}\in \mathcal{F}_k$, $\|X_k^{-1}(\tilde{x} - x^k)\|_2 \le n+2\sqrt{n}$.
\end{lemma}
This lemma is essentially Theorem 5.3.8 of \citep{nesterov2018lectures} in the case when $x$ is on the central path. 
It provides an upper bound of the distance from $x^k$ to any $\tilde{x} $ in the sublevel set of $x^k$. Now we proceed to the proof of \Cref{thm: controlled perturbation for subspace LP}.
\begin{proof}{Proof of \Cref{thm: controlled perturbation for subspace LP}}
To start with, we need to prove that $\hat{x}\in \mathcal{F}_k$.
Suppose that $\hat{x}$ is not in $\mathcal{F}_k$ and let $\tilde{x}$ be the minimizer of $\arg\min_{\bar{x} \in \mathcal{F}_k} (c+p)^\top \bar{x}$. Then due to the convexity of the feasible set,  $\tilde{x}$ must lie on the same constant-objective hyperplane of $x^k$, i.e., $\tilde{x}\in \{x: c^\top x = c^\top x^k\}$. Therefore, it suffices to show that once $p$ satisfies \eqref{eq:sizeofperturbation}, then $\tilde{x}\notin\{x: c^\top x = c^\top x^k\}$, namely, 
\begin{equation}\label{eq thm: controlled perturbation 0}
\min_{\tilde{x}\in \mathcal{F}_k \cap \{x: c^\top x = c^\top x^k\}} \ (c+p)^\top \tilde{x} > \min_{z \in \mathcal{F}_k} \ (c+p)^\top z \ .
\end{equation}
The above inequality \eqref{eq thm: controlled perturbation 0} can be equivalently rewritten as
\begin{equation}\label{eq thm: controlled perturbation 1}
\min_{\tilde{x}\in \mathcal{F}_k \cap \{x: c^\top x = c^\top x^k\}} \ (c+p)^\top (\tilde{x} - x^k) > \min_{z \in \mathcal{F}_k} \ (c+p)^\top  (z - x^k) \ .
\end{equation}
by subtracting $(c+p)^\top x^k$ on both sides of \eqref{eq thm: controlled perturbation 0}. For simplicity of notation, we use LHS$_p$ and RHS$_p$ to denote the left-hand and right-hand sides of \eqref{eq thm: controlled perturbation 1}, respectively. Later we will show $\mathrm{LHS}_p > \mathrm{RHS}_p$ to prove \eqref{eq thm: controlled perturbation 0}. 

Firstly, we derive a lower bound of LHS$_p$. Using $d$ to denote $\tilde{x} - x^k$, we have
\begin{equation}\label{eq thm: controlled perturbation 2}
\mathrm{LHS}_p  = \min_{d\in (\mathcal{F}_k \cap \{x: c^\top x = c^\top x^k\})-x^k} \ (c+p)^\top d = \min_{d\in (\mathcal{F}_k - x^k) \cap \{d: c^\top d = 0\}} \ (c+p)^\top d  = \min_{d\in (\mathcal{F}_k - x^k) \cap \{d: c^\top d = 0\}} \ p^\top d
\end{equation}
where the last equality is due to the linear constraint $c^\top d = 0$.
Furthermore, because of \Cref{lm: analytic center rounding}, $(\mathcal{F}_k - x^k) \cap \{x: c^\top x = 0\}\subseteq \mathcal{F}_k - x^k \subseteq \{d:\|X_k^{-1} d\|_2 \le n+2\sqrt{n}\}$, and then 
\begin{equation}\label{eq thm: controlled perturbation 2-1}
\mathrm{LHS}_p \ge \min_{d:\|X_k^{-1} d\|_2 \le n+2\sqrt{n} } \ p^\top d  = 
\min_{\tilde{d}:\|\tilde{d}\|_2 \le n+2\sqrt{n} } \ (X_k p)^\top \tilde{d}
=  -(n+2\sqrt{n})\|X_k p\|_2 
\end{equation}
in which the first equality holds by substituting $X_k^{-1}d$ with $\tilde{d}$ and the second equality holds because the optimal $\tilde{d}$ is $-\tfrac{(n+2\sqrt{n})\cdot X_k p}{\|X_k p\|_2}$.

Secondly, we derive an upper bound of RHS$_p$. Using $r$ to denote $z - x^k$, we have
\begin{equation}\label{eq thm: controlled perturbation 3}
\mathrm{RHS}_p = \min_{r \in \mathcal{F}_k - x^k} \ (c+p)^\top r = \min_{r \in L' \cap \mathbb{R}^n_+  \cap \{r:c^\top r \le c^\top x^k\}  - x^k} \ (c+p)^\top r  = \min_{r \in L \cap (\mathbb{R}^n_+ - x^k) \cap \{r:c^\top r \le 0\}} \ (c+p)^\top r \ .
\end{equation}
It should be noted that since $\mathbb{R}^n_+ \supseteq \big\{z:\sum_{i=1}^n \big(z_i/x^k_i-1\big)^2\le 1\big\} = \big\{z: \|X_k^{-1}(z - x^k)\|_2\le 1\big\}$, then $\mathbb{R}^n_+ - x^k \supseteq \{r:\|X_k^{-1} r\|_2 \le 1\}$. Therefore, RHS$_p$ can be upper bounded by:
\begin{equation}\label{eq thm: controlled perturbation 3-1}
\mathrm{RHS}_p  \le   \min_{r \in L \cap \{r:\|X_k^{-1} r\|_2 \le 1\} \cap \{r:c^\top r \le 0\}} \ (c+p)^\top r  = c^\top \hat{r} + p^\top \hat{r}
\end{equation}
where $\hat{r}$ denotes the optimal solution $\arg\min_{r \in L \cap \{r:\|X_k^{-1} r\|_2 \le 1\} \cap \{r:c^\top r \le 0\}} \ c^\top r$. The $p^\top \hat{r}$ term has the following upper bound: 
$$
p^\top \hat{r} \le \max_{r \in L \cap \{r:\|X_k^{-1} r\|_2 \le 1\} \cap \{r:c^\top r \le 0\}} \ p^\top r \le \max_{r:\|X_k^{-1} r\|_2 \le 1 } \ p^\top r  = \max_{\tilde{r}:\|\tilde{r}\|_2 \le 1} \ (X_k p)^\top \tilde{r}  = \|X_kp\|_2.
$$
As for the $c^\top \hat{r}$ term in \eqref{eq thm: controlled perturbation 3-1}, it is equal to:
\begin{equation}\label{eq thm: controlled perturbation 3-2}
c^\top \hat{r} = \min_{r \in L \cap \{r:\|X_k^{-1} r\|_2 \le 1\} \cap \{r:c^\top r \le 0\}} \ c^\top r = \min_{r \in L \cap \{r:\|X_k^{-1} r\|_2 \le 1\} } \ c^\top r  = \min_{\tilde{r} \in X_k^{-1} L \cap \{\hat{r}:\|\hat{r}\|_2 \le 1\}} \ (X_k c)^\top \tilde{r}      \ .
\end{equation}
Note that for any $\tilde{r}$ in the linear subspace $X_k^{-1}L$, $(X_k c)^\top \tilde{r} = (P_{X_k^{-1}L}(X_k c))^\top \tilde{r}$ and thus 
$$
\min_{\tilde{r} \in X_k^{-1} L \cap \{\hat{r}:\|\hat{r}\|_2 \le 1\}} \ (X_k c)^\top \tilde{r} = \min_{\tilde{r} \in X_k^{-1} L \cap \{\hat{r}:\|\hat{r}\|_2 \le 1\}} \ (P_{X_k^{-1}L}(X_k c))^\top\tilde{r}  =-\|P_{X_k^{-1}L}(X_k c)\|_2.
$$
Substituting the $c^\top \hat{r}$ term
and the upper bound of $p^\top \hat{r}$ into \eqref{eq thm: controlled perturbation 3-1} yields an upper bound of RHS$_p$: $\mathrm{RHS}_p\le -\|P_{X_k^{-1}L}(X_k c)\|_2 + \|X_kp\|_2$. 

Lastly, when $p$ satisfies \Cref{eq:sizeofperturbation}, we have that
$$
- (n+2\sqrt{n})\cdot\|X_k p\|_2 > -   \| P_{X_k^{-1} \cdot L} (X_k c)\|_2  + \|X_k p\|_2 \ .
$$
Using the lower bound \eqref{eq thm: controlled perturbation 2-1} on $\mathrm{LHS}_p$ and the upper bound we just obtained on $\mathrm{RHS}_p$, we get
$$
\mathrm{LHS}_p \ge -(n+2\sqrt{n})\|X_k p\|_2 > -   \| P_{X_k^{-1} \cdot L} (X_k c)\|_2  + \|X_k p\|_2 \ge \mathrm{RHS}_p \ ,
$$
which proves \Cref{eq thm: controlled perturbation 0} and $c^\top \hat{x} \le c^\top x_k$.\end{proof}


The above results theoretically guarantee the quality of the solutions provided by the smartly perturbed problems. Next, we will show how to integrate optimal face detection and smart perturbations to improve computational efficiency in practice.

\subsection{Perturbation Crossover Method}

In this subsection, we introduce the  ``perturbation crossover method,'' which integrates the smart perturbation and the  candidate optimal face detection.
The overall process is summarized in \Cref{fig:figureofalg}, which focuses on the perturbation crossover applied to the primal problem to obtain a BFS with a primal-dual objective gap smaller than $\varepsilon$.  Note that the corresponding crossover applied to the dual problem is symmetric and will return a dual BFS.

\begin{figure}[htbp]
\centering
\begin{tikzpicture}[-stealth,scale=0.75,transform shape]
\node[end] (2) at (0,0) {Starting Interior-Point solution $(x^k,y^k)$};
\node[decide] (3) [right=of 2] {$\begin{array}{c}
    \text{Feasibility} \\
    \text{Problem?}
\end{array}$};
\draw (2) -- (3);
\node[exec] (5) [right=of 3] {$\begin{array}{c}
    \text{Solve a Randomly} \\
    \text{Perturbed Problem}
\end{array}$};
\draw (3) -- node [above] {Yes} (5);
\node[inout] (e1) [below=of 5] {BFS $x^\star$};
\draw (5) -- (e1);
\node[endn] (6) [below=of e1] {End};
\draw (e1) -- (6);

\node[exec] (7) [below=of 3] {$\begin{array}{c}
       \text{Compute a Perturbation} \\
    \text{that Satisfies \eqref{eq:sizeofperturbation}}
\end{array}$};
\draw (3) -- node [right] {No} (7);
\node[exec] (8) [left=of 7] {$\begin{array}{c}
 \text{Identify a Candidate} \\
    \text{Optimal Face } (\Theta_{p,\gamma}^k, \Theta_{d,\gamma}^k)
\end{array}$};
\draw (7) -- (8);
\node[exec] (9) [left=of 8] {
$\begin{array}{c}
    \text{Solve the Perturbed} \\
    \text{Restricted Problem}
\end{array}$};
\draw (8) -- (9);
\node[decide] (10) [below=of 9] {Infeasible?};
\draw (9) -- (10);
\node[exec] (11) [right=79.3pt of 10] {Decrease $\gamma$};
\draw (10) -- node [above] {Yes} (11);
\draw (11) -- (8);
\node[inout] (12) [below=of 10] {BFS $\hat{x}$};
\draw (10) -- node [right] {No} (12);
\node[decide] (13) [right=of 12] {$c^\top \hat{x}-b^\top y^k < \varepsilon$};
\draw (12) -- (13);
\node[endn] (14) [below=of 13] {End};
\draw (13) -- node [right] {Yes} (14);
\node[exec] (15) [right= of 13] {Reoptimize};
\draw (13) -- node [above] {No} (15);
\node[inout] (e2) [right=of 15] {BFS $x^\star$};
\draw (15) -- (e2);
\node[endn] (16) [right=of e2] {End};
\draw (e2) -- (16);
\end{tikzpicture}
\caption{The flowchart of the perturbation crossover for obtaining an optimal basic feasible solution} \label{fig:figureofalg}
\end{figure}

The first step is to detect whether all feasible solutions are optimal, and we call the LP in this case a feasibility problem. In practice, when $P_{L}(c) = 0$, we say the problem is a feasibility problem.
If the problem is identified as a feasibility problem, according to \Cref{lm: perturbation find BFS,lm: perturbation find BFS dual}, solving any randomly perturbed problem that has an optimal solution almost surely provides an optimal BFS $x^\star$.
If the problem is not a feasibility problem, the crossover proceeds by
computing a smart perturbation that satisfies \eqref{eq:sizeofperturbation} and detecting the candidate optimal face.
Then the BFS is obtained from solving the perturbed restricted problem, which is the LP whose variables not in the optimal face are fixed at their bounds and the objective vector is changed by the smart perturbation. 
If the restricted problem is infeasible, then $\gamma$ in \eqref{eq: relaxed identified optimal face} is decreased and the restricted problem is solved again with the same perturbation but more unfixed variables. 
If it is feasible, the optimal solution $\hat{x}$ of the perturbed restricted problem is a BFS for the original LP,
and the quality of its objective value can be examined by the duality gap between the BFS $\hat{x}$ and the given interior-point dual solution $y^k$, namely $c^\top \hat{x} - b^\top y^k$.
If the duality gap is larger than the tolerance $\varepsilon$, 
the reoptimization phase follows by warm starting the simplex method from $\hat{x}$. 

In practice, the original LP may not be in standard form,  so we transform the problem into standard form before computing the smart perturbation. Subsequently, after identifying the optimal face, we apply the restrictions and the perturbation to the original LP and let an LP solver solve it. In this framework, solving one linear system is required when checking whether the LP is a feasibility problem and computing the perturbation. However, since the interior-point method solves a linear system in each iteration (and may have already stored the needed matrix factorization), the cost of solving the additional linear system incurred by our approach is negligible.

It should be mentioned that \cite{mehrotra1991finding} also proposed an innovative approach to identify corner solutions via perturbations. Mehrotra's method specifically focused on primal perturbations for \Cref{pro: General LP}. While their method provided valuable insights, we provide a theoretical guarantee for managing the magnitude of the perturbation, ensuring the high quality of the perturbed problem's solutions. In practice, our method additionally restricts the variables inside a candidate optimal face, significantly reducing the computation time of the perturbed problem.
\section{Numerical Experiments}\label{sec: Numerical Experiments}

In \Cref{subsec: network experiments}, we first evaluate our network crossover method~(\Cref{sec: Network Crossover Method}) on optimal transport and MCF problems with three datasets: 1) MNIST image data\footnote{\url{https://www.tensorflow.org/datasets/catalog/mnist}. Accessed: 20 May 2023. This access date applies to all datasets mentioned below.}, 2) Hans Mittelmann's benchmark problems\footnote{\url{https://plato.asu.edu/ftp/network.html}}, and 3) GOTO (Grid on Torus) instances\footnote{\url{https://lemon.cs.elte.hu/trac/lemon/wiki/MinCostFlowData}}. 
We follow by evaluating the perturbation crossover method~(\Cref{sec: Perturbation Crossover Method}) on problems from Hans Mittelmann's benchmarks for optimization software\footnote{\url{https://plato.asu.edu/ftp/lpopt.html}} in \Cref{subsec: perturbation experiments}. Details of how these data sets are collected and used can be found in the Appendices.
We run and compare our methods with the typical commercial solvers: \textsc{Gurobi} (9.0.3),
\textsc{Cplex} (Studio 12.9), and \textsc{Mosek} (9.3.14).
The experiments are implemented using Python 3.7 and run on a macOS system~(Ventura 13.0) with a hardware configuration of 64 GB RAM and Apple M1 Max chip.

All our codes and data are available from the IJOC GitHub software repository \citep{SmartCrossover}.

\subsection{Network Crossover Method} \label{subsec: network experiments}

In \Cref{sec: Network Crossover Method}, we propose \textsc{Tree-BI} and \textsc{Col-BI} algorithms for the basis identification step of the network crossover method. 
For ease of reference, we denote \textsc{TNET} and \textsc{CNET} as the network crossover method when employing these two basis identification methods and the corresponding reoptimization phase. First, we show the experimental results on optimal transport problems. 

\subsubsection{Optimal Transport.} \label{subsubsec: OT experiments}
We conduct experiments on the optimal transport problem because it has extensively-studied first-order methods available, making it an ideal testing ground for assessing not only the crossover's efficacy but also the power of combining crossover and first-order methods.

\paragraph{Experiment Configurations.}
For \textsc{Col-BI}, we set $\theta_k$ in \Cref{eq: generate series of subproblem}, the value controlling the number of entering columns, as $2^k$. 
We set the regularization coefficient in the \textsc{Sinkhorn} algorithm ($\lambda$ in \citep{cuturi2013sinkhorn}) to 10.
In the \textsc{Col-BI} in \textsc{CNET} and the reoptimization stage in \textsc{TNET} \& \textsc{CNET}, we repeatedly call the simplex method from the corresponding commercial solver we use for comparison. For simplicity of notations, we use \textsc{TNET-grb} and \textsc{CNET-grb}, to denote \textsc{TNET} and \textsc{CNET} with \textsc{Gurobi}'s simplex method. Similarly, we use \textsc{TNET-cpl} and \textsc{CNET-cpl} to denote the same scheme with \textsc{Cplex}'s network simplex method. 
For fairness of comparison, we set the pricing strategy to use the steepest edge, to exclude the impact of heuristic pricing. 
We also turn off the $\textit{presolve}$ because the presolve only removes one constraint for each optimal transport problem and does not change the solving time much, while it takes a very long time for large instances.

\begin{figure}[htb]
    \centering
    \includegraphics[width=0.7\textwidth]{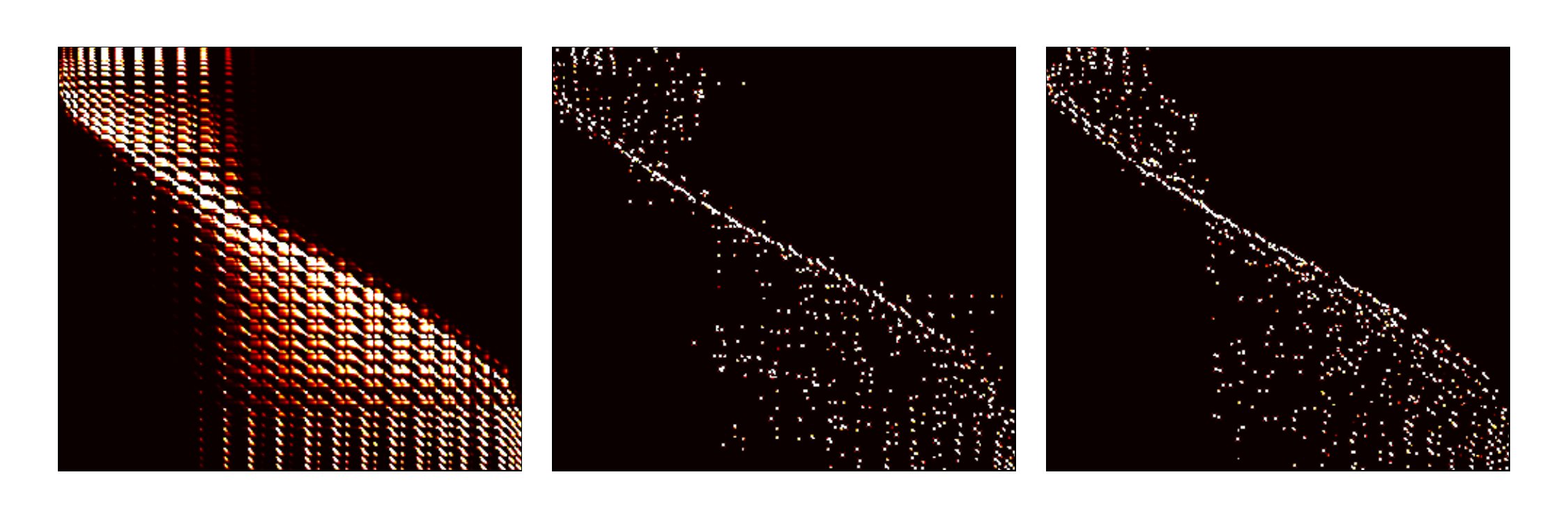}
    \vspace{-5pt}
    \caption{Transport plans of an optimal transport problem. The left, the middle, and the right subfigures are the transportation plans with respect to the initial interior-point solution, the basic solution from \textsc{Tree-BI}, and the optimal BFS after reoptimization. The solution in the left subfigure is generated by \textsc{Gurobi}'s barrier algorithm with the barrier convergence tolerance set to 0.01. In these figures, the brightness of each point at $(i,j)$ represents the amount of grayscale transported from the $i$-th pixel to the $j$-th pixel.}
    \label{fig: TransportPlan}
\end{figure}

First of all, we show that the basic solution found via \textsc{Tree-BI} presents much higher sparsity while maintaining similarity to the solution of the starting method.
For an optimal transport problem from two randomly selected images, \Cref{fig: TransportPlan} shows the weight of the nonzeros in the initial interior-point solution, the basic solution found by \textsc{Tree-BI} and the output after reoptimization.  One can observe that the BFS found by \textsc{Tree-BI} can promote higher sparsity while closely approximating the key information of the interior-point solution and the optimal BFS.

Next, our main experiments will examine the two following questions:
\begin{itemize}
\item[a)] The crossover procedure in isolation. 
We use the solutions from \textsc{Gurobi}'s barrier method, with different convergence tolerances, to approximate a solution from the starting method. Based on this solution, we compare \textsc{TNET}, \textsc{CNET}, and \textsc{Gurobi}'s internal crossover stage. 
\item[b)] The overall running time of obtaining an optimal BFS (the time of the crossover procedure coupled with the preceding interior-point or first-order method). The methods in the comparison include \textsc{Gurobi}'s default method, \textsc{Cplex}'s network simplex method, and the \textsc{Sinkhorn} algorithm augmented with \textsc{TNET} or \textsc{CNET}. 
\end{itemize}

First, we examine (a) by comparing different crossover methods. 
\Cref{tab: comparecrosstime} compares crossover methods on optimal transport problems of different scales, by starting from the same interior-point solution returned by \textsc{Gurobi}. 
Both \textsc{TNET} and \textsc{CNET} show improvements compared to \textsc{Gurobi}'s crossover in larger instances.
Interestingly, \textsc{TNET} exhibits a larger advantage when the given initial interior-point solution has higher accuracy.

\begin{table}[htbp]
  \centering
   \caption{Crossover procedure comparison of \textsc{Gurobi}'s barrier part (grbBarrier), \textsc{Gurobi}'s  crossover part (grbCrossover), \textsc{CNET-grb} and \textsc{TNET-grb} on optimal transport problems generated from MNIST dataset. For each scalarization factor $\alpha$, we test on 10 independent trials and report the geometric average running time (in seconds). The barrier convergence tolerance of the interior-point solution is set to be 1e-8 (the left subtable) and 1e-2 (the right subtable), respectively.}
   \vspace{6pt}
  \setlength{\tabcolsep}{2mm}{
  \scalebox{0.7}{
    \begin{tabular}{rrrrr}
    \toprule
      $\alpha$  &  grbBarrier &  grbCrossover &  \textsc{CNET-grb}  &  \textsc{TNET-grb} \\
    \midrule
      1 &       0.12 &       0.02  &      0.02     &      0.02 \\
      2 &       2.05  &       0.11  &     \textbf{0.14}     &      0.15 \\
      3 &      11.23  &       0.66  &     \textbf{0.64}     &      0.73  \\
      4 &      50.60  &      8.74  &     \textbf{3.78}     &      3.89   \\
      5 &      79.45  &      29.55  &      5.70     &     \textbf{4.82}   \\
      6 &     165.71  &     98.24  &     13.53     &   \textbf{11.00}    \\
      7 &     813.71  &     809.67 &     60.75     &   \textbf{53.55}    \\
    \bottomrule
    \end{tabular}
    \quad
    \begin{tabular}{rrrrr}
    \toprule
      $\alpha$  &  grbBarrier &  grbCrossover &  \textsc{CNET-grb}  &  \textsc{TNET-grb} \\
    \midrule
       1 &        0.10  &       0.02  &      0.02  &      0.02 \\
       2 &        1.70  &       0.22 &      \textbf{0.14}  &       0.15  \\
       3 &       9.19   &       1.45  & \textbf{0.65}  &      0.74  \\
       4 &       40.69  &      22.07  &  \textbf{3.76}  &      4.12   \\
       5 &       59.83  &      37.47  &      5.38  & \textbf{5.04}   \\
       6 &      126.31  &     125.26 &     14.70  &  \textbf{11.69}  \\
       7 &      570.56  &     1015.20  &    76.75  &  \textbf{72.85}    \\
    \bottomrule
    \end{tabular}
    }}
  \label{tab: comparecrosstime}
\end{table}

Second, we examine (b) by comparing the overall running time of different methods for obtaining an optimal BFS.  
As \Cref{tab: compareTtime_mnist}  shows, the Sinkhorn+\textsc{CNET-cpl} is slightly better when compared with Sinkhorn+\textsc{TNET-cpl}, and often faster than the network simplex method, particularly in larger-scale scenarios. 
These advantages also show the potential of combining other first-order methods and network simplex methods, compared with directly running the network simplex method alone. We use the nextwork simplex method implementation of \textsc{Cplex} for its warm-start capability, which facilitates our experimental setup. In practice, the combination with other more efficient network simplex method implementations, such as those in the LEMON library\footnote{\url{https://lemon.cs.elte.hu/trac/lemon}} or the latest version of \textsc{Gurobi}, may further enhance the computational speed of all these methods.

\begin{table}[htbp]
  \centering
   \caption{Total running time comparison of \textsc{Gurobi}'s default method (grbDefault), \textsc{Cplex}'s network simplex method (cplNetSimplex), \textsc{Sinkhorn}+\textsc{CNET-cpl}, and \textsc{Sinkhorn}+\textsc{TNET-cpl} on optimal transport problems generated from MNIST. 
   For each scalarization factor $\alpha$, we test on 10 independent trials and report the geometric average running time (in seconds) and iteration number, which includes all simplex iterations applied, and the push iterations to find a feasible tree BFS (only for \textsc{TNET}). }
   \vspace{8pt}
  \setlength{\tabcolsep}{3.5mm}{
  \scalebox{0.75}{
    \begin{tabular}{rrrrrrrr}
    \toprule
    \multicolumn{1}{c}{\multirow{2}{*}{$\alpha$}} & \multicolumn{1}{c}{\multirow{2}{*}{grbDefault}} & \multicolumn{2}{c}{cplNetSimplex} & \multicolumn{2}{c}{Sinkhorn+\textsc{CNET-cpl}} & \multicolumn{2}{c}{Sinkhorn+\textsc{TNET-cpl}} \\
    \cmidrule(r){3-4} \cmidrule(r){5-6} \cmidrule(r){7-8} 
    \multicolumn{1}{c}{} & \multicolumn{1}{c}{} & \multicolumn{1}{c}{time} &  \multicolumn{1}{c}{iterations} &  \multicolumn{1}{c}{time} &  \multicolumn{1}{c}{iterations} & \multicolumn{1}{c}{time} &  \multicolumn{1}{c}{iterations} \\
    \midrule
      1 &  0.02       &     0.01       &  1.6e+03  &      0.02       & 1.5e+03 & 0.02   & 6.0e+02 \\
      2 &  0.24       &  \textbf{0.03} &  1.2e+04  &      0.10       & 1.0e+04 & 0.17   & 5.1e+03 \\
      3 &  2.43       &  \textbf{0.34} &  6.1e+04  &      0.50       & 3.4e+04 & 0.84   & 1.9e+04 \\
      4 &  12.60      &     2.87       &  2.0e+05  &  \textbf{2.18}  & 9.4e+04 & 3.38   & 5.5e+04 \\
      5 &  40.47      &     4.84       &  3.0e+05  &  \textbf{3.74}  & 1.2e+05 & 5.78   & 8.9e+05 \\
      6 &  245.34     &    12.84       &  6.3e+05  &  \textbf{8.54}  & 2.6e+05 & 12.50  & 1.6e+05 \\
      7 &  1521.80    &    62.10       &  1.6e+06  &  \textbf{37.14} & 6.1e+05 & 54.12  & 3.9e+05 \\
  \bottomrule
    \end{tabular}}}
  \label{tab: compareTtime_mnist}
\end{table}%

\subsubsection{Minimum Cost Flow.} \label{subsubsec: MCF experiments}
Next, we extend experiments to more general MCF problems. These experiments still use the solutions from early-stopping \textsc{Gurobi}'s barrier method (to approximate the solution of the starting method) and  compare our crossover methods with \textsc{Gurobi}'s crossover.
We conclude by noting that the potential value of our approach lies particularly in the context of solving large-scale MCF problems in which the solution given by the starting method is of low accuracy,  such as  outputs of first-order methods in large-scale LP tasks.

We conducted our experiments using benchmark problems from two sources: 1) Hans Mittelmann's large network-LP benchmark\footnote{\url{http://plato.asu.edu/ftp/lptestset/network/}}, and 2) networks generated by GOTO\footnote{\url{http://lime.cs.elte.hu/~kpeter/data/mcf/goto/}.}. 
The experiment configurations are the same as the experiments on optimal transport problems in \Cref{subsubsec: OT experiments}.

\begin{figure}[htb]
    \centering
    \includegraphics[width=0.9\textwidth]{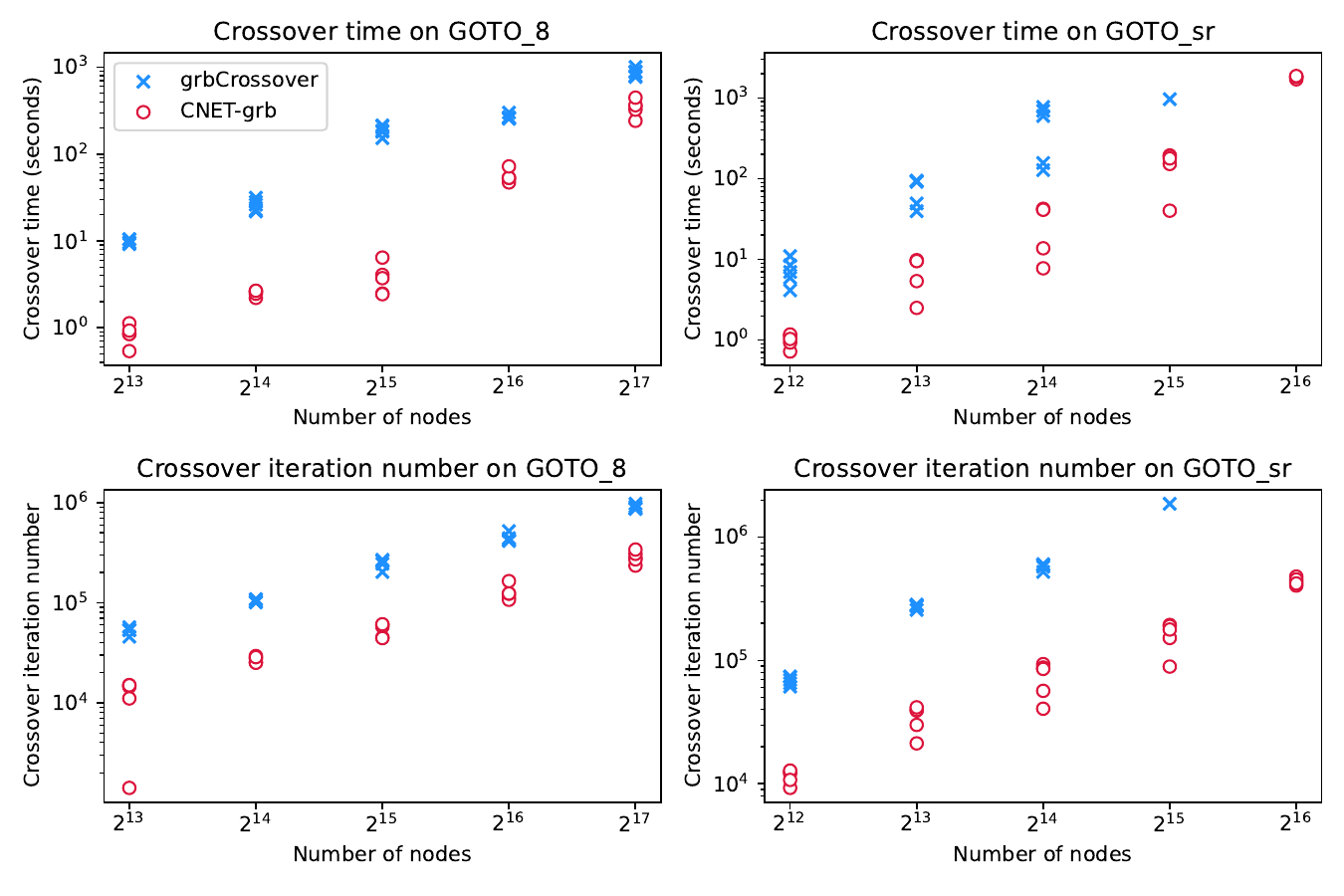}
    \vspace{-5pt}
    \caption{Comparison of \textsc{Gurobi}'s crossover and \textsc{CNET-grb} on GOTO\_8 and GOTO\_sr with different scales. The barrier convergence tolerance of the interior point solution is set to be 1e-2. The iterations of \textsc{Gurobi} are all the iterations in the crossover (push phase and reoptimization) and the iterations of \textsc{CNET-grb} are all iterations of \textsc{Col-BI} and reoptimization. For each problem scale, we report the results of 5 different generated problems. The time limit is set as 4000 seconds.}
    \label{fig: goto}
\end{figure}

\Cref{tab: mcf crossover} shows that on the large network-LP benchmark dataset, when the given solution is of low accuracy, \textsc{CNET} is much better than \textsc{Gurobi}'s crossover in both overall running time and iteration number. When the given solution is of high accuracy, \textsc{Gurobi}'s crossover has advantages, but both crossover methods are very fast (no more than 10 seconds) in this case.  
\Cref{fig: goto} compares \textsc{CNET} with \textsc{Gurobi}'s crossover with the problem scale varying from small to large, in which the largest problems (GOTO\_sr with $2^{16}$ nodes) contain over 10 million variables. See the Appendix for the complete tables of the experiments in \Cref{tab: mcf crossover} and \Cref{fig: goto}.
These results demonstrate the consistent advantages of our crossover methods over \textsc{Gurobi}'s crossover when the given solution from the starting method is of low accuracy (for example, solution of first-order methods). 
In general, \textsc{Gurobi}'s crossover heavily relies on high-accuracy input solutions while \textsc{CNET} is more stable. This demonstrates that the graph information used in \textsc{CNET} allows for a more robust crossover, which is particularly advantageous when starting from  solutions obtained with first-order methods. 

\begin{table}[htbp]  
    \centering
    \caption{
    The comparison of \textsc{Gurobi}'s crossover (grbCrossover) and \textsc{CNET-grb} on large network-LP benchmark dataset. The primal-dual gap of the interior-point solution is set to be 1e-2 and 1e-8. The ``GeoAvg Time'' and ``GeoAvg Iterations'' denote the geometric average of running time and iteration number for all problems.
    }
    \vspace{4pt}
    \setlength{\tabcolsep}{4.8mm}{
    \scalebox{0.75}{
        \begin{tabular}{cccccc}
            \toprule
            \multicolumn{1}{c}{\multirow{2}{*}{gap}} & \multicolumn{1}{c}{\multirow{2}{*}{grbBarrier}} & \multicolumn{2}{c}{grbCrossover} & \multicolumn{2}{c}{\textsc{CNET-grb}} \\
            \cmidrule(r){3-4} \cmidrule(r){5-6} 
            \multicolumn{1}{c}{} & \multicolumn{1}{c}{} &  \multicolumn{1}{c}{GeoAvg Time} &  \multicolumn{1}{c}{GeoAvg Iterations} &  \multicolumn{1}{c}{GeoAvg Time} &  \multicolumn{1}{c}{GeoAvg Iterations} \\
            \midrule
            1e-2 & 395.6  &   68.0   & 9.8e+05   &   8.40  &  4.8e+04 \\
            1e-8 & 581.1 &  0.52 & 5.5e+03   & 4.92 & 3.9e+04  \\ 
            \hline
        \end{tabular}
        }}
    \label{tab: mcf crossover}
\end{table}

\subsection{Perturbation Crossover Method}\label{subsec: perturbation experiments}

In this subsection, we describe the implementation of the perturbation crossover method, in which the perturbed problems are solved by commercial solvers. We then compare it with the solvers' internal crossover. 
The experiment setup for the perturbation crossover on one LP is as follows:
\begin{itemize}
    \item[\textit{Step 1}] 
    \textit{Preprocessing}. For the given LP, we apply \textsc{Gurobi}'s presolve function and obtain a presolved model. In LP solving, the presolve function reduces the problem size by removing redundancies and eliminating variables. It changes the model and advances the realistic solving substantially. Thus, by starting with presolve, we ensure that our experimental framework reflects a realistic and practical optimization setting. 
    \item[\textit{Step 2}] \textit{Barrier method of the solver}. We first employ the solver's barrier method to solve the problem, and the resulting interior-point solution will be used for the perturbation crossover method. 
    We follow the solver's default tolerance for all problems.
    \item[\textit{Step 3}]
    \textit{Perturbation crossover.} We adhere to the procedure outlined in the flowchart presented in \Cref{fig:figureofalg}. The initial parameter $\gamma$, for estimating the candidate optimal face, is set to $\expnumber{1}{-3}$. The perturbation vector $p$ on the objective vector, inspired by  \Cref{thm: controlled perturbation for subspace LP}, is calculated on the standard form reformulation and applied more aggressively in practice. For the new artificial variables in the standard form, $p_i$ is zero. For other indices $i$, the perturbation is defined as:
    \begin{equation} \label{eq: perturb in experiment}
        p_i := \frac{\xi_i}{\| \xi \|_2} \cdot \frac{\|P_{X^{-1}\cdot L}(Xc)\|_2}{0.01 \cdot n \cdot \max\{10^{-6},x_i\}}, \quad \text{for all } i\in [n].
    \end{equation}
    Here, $\xi$ is a random vector from a uniform distribution between $0.9$ and $1$, $x = (x_i)_{i \in [n]}$ represents the interior-point primal solution, and $P_{X^{-1}\cdot L}(Xc)$ is as defined in \Cref{eq: projector def}. 
    The construction of \Cref{eq: perturb in experiment} follows the insight provided by \Cref{eq:sizeofperturbation} of \Cref{thm: controlled perturbation for subspace LP}. The division by $0.01$ is to intensify the perturbation, a strategy that has been empirically proven more efficient in our tests. Taking the maximum with $10^{-6}$ is aimed at mitigating numerical issues with small values of $x_i$.
    In cases where the restricted problem is infeasible,\footnote{This is not common. See Appendix \Cref{app: perturb experiment} for more information.} $\gamma$ is reduced by $\gamma \leftarrow \gamma \cdot \expnumber{1}{-5}$ repeatedly until the restricted problem becomes feasible. The presolve function is activated when solving the perturbed problem.
   
    \item[\textit{Step 4}] \textit{Reoptimization}. When the relative objective gap between $f_p$ (the objective value of the solution of the perturbed restricted problem) and $f_d$ (the objective value of the given dual interior-point solution), defined as $\frac{\lvert f_p - f_d\rvert}{\lvert f_p\rvert + \lvert f_d\rvert + 1}$, is smaller than the default tolerance $\varepsilon:=\expnumber{1}{-8}$, then the BFS obtained is considered sufficiently optimal.
    This tolerance is similar to those of most commercial solvers, like \texttt{BarConvTol} in \textsc{Gurobi} and \texttt{intpntTolRelGap} in \textsc{Mosek}.
    If the gap exceeds the tolerance, the reoptimization follows by warm starting the simplex method from the BFS. 
\end{itemize}


We measure the running time of the perturbation method by tracking the total elapsed time taken to solve the perturbed model and the time for the following reoptimization phase. 
The time consumed in determining the feasibility of the problem and in computing the perturbations would be negligible if these processes are integrated within the solvers internally. Consequently, these time components are excluded from the reported running time of the perturbation crossover method. 

The first experiment is carried out on the LP relaxation of five max cut packing problems from \cite{bergner2014branch}. These problems are also collected in the MIPLIB 2017 dataset \citep{gleixner2021miplib}. For the original mixed integer program, a BFS of the LP relaxation is critical for the subsequent branches and cuts. \Cref{tab:maxcutpacking} shows the runtime for \textsc{Gurobi}'s barrier method and crossover, and our perturbation crossover, along with the basic information of the problems.
As shown in \Cref{tab:maxcutpacking}, the five problems vary not only in scale but also in the dimensions of their optimal faces. While \textsc{Gurobi} is consistently fast in the barrier method part, it tends to consume a considerable amount of time during the crossover phase, especially when the (estimated) dimension of the optimal face is large.
When this dimension is large, many linear constraints are inactive for the given interior-point solution, and thus the traditional crossover spends more time identifying an optimal basis. In contrast, our perturbation crossover significantly outperforms \textsc{Gurobi}'s crossover in these cases.

\begin{table}[htbp]
    \centering
    \caption{The experimental results of $5$ max cut packing problems. ``Est. Dim.'' denotes the estimated dimension of the optimal face. ``grbBarrier'', ``grbCrossover'', and ``perturbCrossover'' show running time (in seconds) of \textsc{Gurobi}'s barrier methods, \textsc{Gurobi}'s crossover methods, and the perturbation crossover method respectively.} \label{tab:maxcutpacking}
    \begin{threeparttable}
    \setlength{\tabcolsep}{3mm}{
    \scalebox{0.75}{
        \begin{tabular}{lcccccc}
        \toprule
        \multicolumn{1}{c}{\multirow{2}{*}{Problem}}  &  \multicolumn{3}{c}{Instance Statistics} &
        \multicolumn{1}{c}{\multirow{2}{*}{grbBarrier}} & \multicolumn{2}{c}{Crossover Time}  \\
        \cmidrule(r){2-4} \cmidrule(r){6-7}
        \multicolumn{1}{c}{} & \multicolumn{1}{c}{Constraints} & \multicolumn{1}{c}{Variables} & \multicolumn{1}{c}{Est. Dim.\tnote{1}} & \multicolumn{1}{c}{} & \multicolumn{1}{c}{grbCrossover} & \multicolumn{1}{c}{perturbCrossover} \\ 
        \midrule
         graph20-20-1rand  &   5587  &  2183  &  2035 &     0.01 &             0.05 &      \textbf{0.04} \\
         graph20-80-1rand  &   55107 &  16263 &  15912&      0.05 &             2.42 &      \textbf{1.11} \\
         graph40-20-1rand  &   99067 &  31243 & 30772 &       0.09 &            15.84 &      \textbf{8.33} \\
         graph40-40-1rand  &  360900 & 102600 & 101700&      0.41 &           323.41 &     \textbf{50.79} \\
         graph40-80-1rand  &  1050112& 283648 & 282112&      1.4  &           $>$10000  &    \textbf{872.07} \\
        \bottomrule
        \end{tabular}
    }}
    \begin{tablenotes}
        \scriptsize
         \item[1] For the ``optimal face'' determined by \Cref{eq: relaxed identified optimal face} with $\gamma=1$, the estimated dimension is calculated by $n$ minus the number of \\ linear equality constraints. 
     \end{tablenotes}
    \end{threeparttable}
\end{table}

We then compare the perturbation crossover method with solvers' internal crossover on Hans Mittelmann's \textit{LPopt} Benchmark problems.\footnote{\url{http://plato.asu.edu/ftp/lpopt.html}} 
To validate the effectiveness of the perturbation crossover across different solvers, we also conducted the same experiments using another widely used commercial solver, \textsc{Mosek} (9.2).
Different from \textsc{Gurobi}, \textsc{Mosek} applies interior-point method on the homogeneous self-dual model \citep{andersen2009homogeneous} so the iterates follow a different central path. For this reason, \Cref{thm: controlled perturbation for subspace LP} cannot directly apply on \textsc{Mosek}.
See \Cref{tab: grb perturb complete} and \Cref{tab: msk perturb complete} of the Appendix for the complete and comprehensive results of the perturbation crossover based on  \textsc{Gurobi} and \textsc{Mosek}.

\begin{figure}[h]
    \centering
    \includegraphics[width=0.9\textwidth]{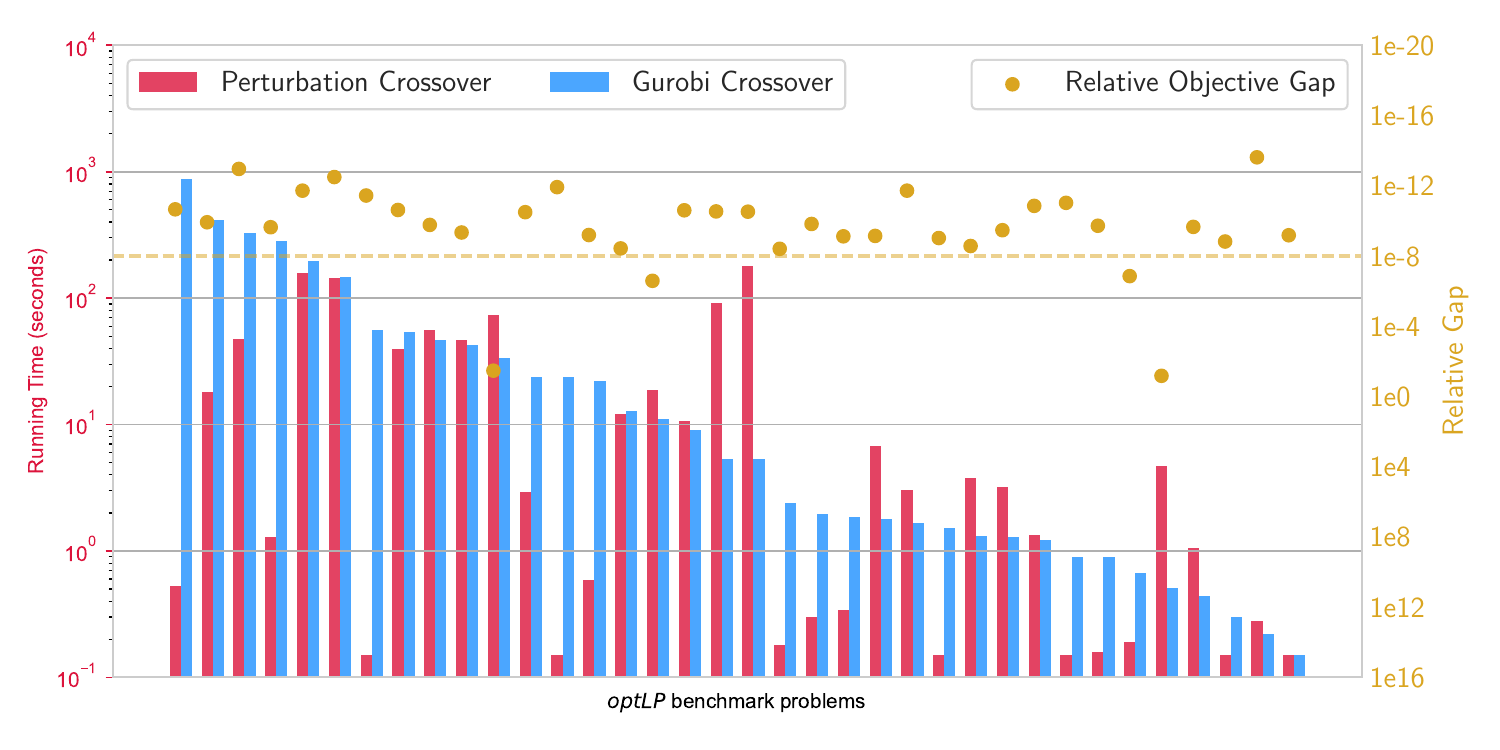}
    \label{fig: grb perturb}
    \caption{Results of the perturbation crossover and \textsc{Gurobi}'s crossover. The x-axis denotes the benchmark problems, arranged according to the crossover time of \textsc{Gurobi}. The red and blue bars (left y-axis) represent the perturbed model's solving time and the solver's original crossover time. The orange dots (right y-axis) represents the relative objective gap between the objective value of the solution of the perturbed restricted problem, and the objective value of the given dual interior-point solution.
    }

    \label{fig: perturb results}
\end{figure}

\Cref{fig: perturb results} presents the time for solving the perturbed restricted problem and \textsc{Gurobi}'s internal crossover, as well as the corresponding relative objective gap of the BFS obtained. One can observe that if the goal is to obtain a BFS within the default tolerance (achieving a relative objective gap smaller than $10^{-8}$), the  BFS obtained is generally accurate enough, with exception in only 4 instances. Furthermore, the time taken to acquire such a BFS is frequently considerably less than that required by the solver's own crossover method.
It should be mentioned that there is a notable variation in the time for solving the perturbed restricted problem and the time for \textsc{Gurobi}' internal crossover, which implies in practice if the computing resource is unlimited (which is often seen in practical applications), concurrently running the two crossover methods yields the shortest running time. 

\Cref{tab: perturb statistics} provides the statistics of the total running time of the perturbation crossover, which includes the time spent on the reoptimization phase when the relative objective gap is not within the tolerance.  Specifically, we exclude problems that fail in the barrier algorithm. This includes cases in which \textsc{Gurobi} automatically switches to the homogeneous barrier algorithm or the barrier algorithm exceeds one hour. 

\begin{table}[h!]
  \centering
    \caption{Performance statistics of the perturbation crossover method using \textsc{Gurobi} and \textsc{Mosek} on finding a primal BFS within the default tolerance in 36 instances. ``Faster'' and ``Slower'' count how many instances get faster when running in our perturbation crossover methods. 
    ``Best per Instance'' denotes the minimum of the solver's original crossover running time or the perturbation crossover method's running time, to approximate concurrently running the two methods.}
  \begin{threeparttable}
    \setlength{\tabcolsep}{3.5mm}{
    \scalebox{0.85}{ 
      \begin{tabular}{lcccccc}
            \toprule
            \multicolumn{1}{c}{\multirow{2}{*}{Solver}} & \multicolumn{1}{c}{\multirow{2}{*}{Faster}} & \multicolumn{1}{c}{\multirow{2}{*}{Slower}} & \multicolumn{3}{c}{GeoAvg Running Time} \\
            \cmidrule{4-6}
            \multicolumn{1}{c}{} & \multicolumn{1}{c}{} & \multicolumn{1}{c}{} & \multicolumn{1}{c}{Solver's Crossover} & \multicolumn{1}{c}{Perturbation Crossover} & \multicolumn{1}{c}{Best per Instance} \\
            \hline
                \textsc{Gurobi}  & 19 & 17 & 6.83 & 3.44  & \textbf{1.53} \\
                \textsc{Mosek}  & 15 & 21  &  5.16 & 10.00 & \textbf{2.03} \\
            \bottomrule
        \end{tabular}
        }}
  \end{threeparttable}
  \label{tab: perturb statistics}
\end{table} 



Although the perturbation crossover method is not consistently better than the traditional crossover implemented in solvers, it presents advantages on some problems, especially those used to be considered challenging, such as the 4 instances that \textsc{Gurobi} spends the most time on (see \Cref{fig: perturb results}). Since \textsc{Mosek} follows a different central path, it makes sense that in \textsc{Mosek} our perturbation crossover does not work as well as the internal crossover, which in turn validates \Cref{thm: controlled perturbation for subspace LP}.
Note that when given sufficient computing resources, running both the traditional crossover and our perturbation crossover methods concurrently can greatly improve computational efficiency and robustness. This is evidenced by the ``Best per Instance'' column, which exhibits a remarkable speedup — approximately 4 to 5 times faster than \textsc{Gurobi}'s internal crossover alone. Moreover, this technique also applies  to \textsc{Mosek}, in which the ``Best per Instance'' is 2 to 3 times faster than the internal crossover, despite the perturbation crossover itself not outperforming the internal crossover of \textsc{Mosek}. 
It is important to note that different solvers employ varied parameter setups and algorithms. For example, the barrier methods used by different solvers may follow distinct central paths because of different internal reformulations. As a result, the results in \textsc{Gurobi} and \textsc{Mosek} are different. Despite these discrepancies, the experimental results obtained thus far strongly support the broad applicability and generalizability of the perturbation crossover approach.

Through the experiments comparing our new perturbation crossover method with the internal crossover methods of commercial solvers, we have demonstrated the distinct advantages of our approach. Our crossover method showcases its capability to navigate issues that have historically challenged traditional crossover methods. Moreover, our results reveal the potential for significant speed enhancements when running both the traditional and perturbation crossover methods concurrently. This indicates a promising direction for future optimization methods and tools.


    \ACKNOWLEDGMENT{%
    This research is partially supported by the National Natural Science Foundation of China (NSFC) [Grant NSFC-72150001, 72225009, 72394360, 72394365].
    We thank Fangkun Qiu for some valuable experiments in the preliminary version of this paper. We thank Jiayi Guo and the developers of the solver COPT, especially Qi Huangfu, for fruitful discussions. We also thank the three anonymous referees for their helpful comments and suggestions, which largely improved the quality of the manuscript. Chengwenjian Wang and Zikai Xiong are corresponding authors, and inquiries may be directed to them.
    }

    \bibliographystyle{informs2014}
    \bibliography{reference} 
    

\begin{APPENDICES}

%
\section{Proof of \Cref{lm: perturbation find BFS}}\label{app:proof}
We present a proof of \Cref{lm: perturbation find BFS} for completeness. 
The proof of \Cref{lm: perturbation find BFS dual} is analogous with the change only in the formulation of the objective vector, so it is omitted. 
\begin{proof}{Proof of \Cref{lm: perturbation find BFS}}
    Let $\mathcal{F}$ denote the feasible set of \eqref{pro: general perturbed primal problem}. It is known from the definition that the extreme points and extreme rays of the optimal face are also the extreme points and extreme rays of the feasible set $\mathcal{F}$. Therefore, if there exist multiple optimal solutions, then there must exist either multiple optimal extreme points of $\mathcal{F}$ or an extreme ray orthogonal with $c+p$. Since $c+p$ has to be orthogonal with the straight line spanned by any two optimal extreme points or the extreme ray, the set of all these $c+p$ is a subset of measure zero in the full space, because they lie in a subspace of a strictly smaller dimension. Overall, since there are only finitely many spanned lines of two extreme points and extreme rays of $\mathcal{F}$, the set of all $c+p$ orthogonal to any of them is still of zero measure. Any $c+p$ such that the perturbed problem has multiple optimal solutions belongs to the above set, so the set of all $p$ so that the perturbed problem has multiple optional solutions is of zero measure. 
\end{proof}

\section{Complete Numerical Results for Network Problems} \label{app: network experiments}

\paragraph{MNIST Dataset.}
We experiment with the MNIST dataset for the optimal transport problem. This dataset contains 60,000 images of handwritten digits ($28\times28$ pixels for each), from which an optimal transport problem is generated by first randomly selecting two images and normalizing the sum of their non-zero grayscale values to 1. And then the optimal transport problem is to find an optimal transport plan from one image (as a discrete distribution with $28\times28$ supports) to the other, with the transport cost defined as the Manhattan distance between the corresponding pixels. To illustrate, transporting the grayscale value from a pixel at location $(7, 9)$ to a pixel at index $(13, 10)$ incurs a cost of $\| (7,9) - (13,10) \|_1 = 7$ per unit.
Furthermore, to expand the problem scale, we artificially increase the image size by a factor of $\alpha$ ($\alpha = 1, 2, \dots, 7$) and create images sized at $28\alpha \times 28 \alpha$. The creation is to split each original pixel into $\alpha \times \alpha$ pixels with the same grayscale and normalize as before. The problems of the largest generated group ($\alpha=7$) have on average 60 million variables.

\paragraph{GOTO Problems.} 
GOTO (Grid on Torus) is one of the hardest problems sets on MCF. We generate our GOTO problems in the same way as \citet{kovacs2015minimum} do. 
Two types of instances are considered: 
1) \textit{GOTO\_8}: sparse network: $m=8n$; 2) \textit{GOTO\_sr}: dense network: $m \approx n\sqrt{n}$, where $m, n$ are the number of arcs and nodes in the network respectively. Note that GOTO problems have been included in Han's Mittelmann's Network LP benchmark problems. Since we do experiments on GOTO separately, we ignore all GOTO problems when considering benchmark problems in \Cref{tab: mcf crossover} and \Cref{tab: complete mcf benchmark results}.

\begin{table}[htbp]
  \centering
   \caption{Complete results of \Cref{tab: mcf crossover} with a starting interior-point solution of primal-dual gap $1e-2$.}
   \vspace{5pt}
  \setlength{\tabcolsep}{4mm}{
  \scalebox{0.8}{
    \begin{tabular}{rrrrrrrr}
    \toprule
    \multicolumn{1}{c}{\multirow{2}{*}{Problem}} & \multicolumn{1}{c}{\multirow{2}{*}{Node}} & \multicolumn{1}{c}{\multirow{2}{*}{Arc}} &  \multicolumn{1}{c}{\multirow{2}{*}{grbBarrier}} & \multicolumn{2}{c}{grbCrossover} & \multicolumn{2}{c}{CNET-grb} \\
    \cmidrule(r){5-6} \cmidrule(r){7-8} 
    \multicolumn{1}{c}{} & \multicolumn{1}{c}{} & \multicolumn{1}{c}{} & \multicolumn{1}{c}{} & \multicolumn{1}{c}{time} &  \multicolumn{1}{c}{iterations} &  \multicolumn{1}{c}{time} &  \multicolumn{1}{c}{iterations} \\
    \midrule
        16\_n14  &  16381 & 261873  &  31.09  &  17.1  & 1.8e+05 &  6.01 &  40547 \\
        i\_n13   &  8181  & 739733  &  24.98  & 100.8  & 2.7e+05 &  7.84 &  38642 \\
        lo10     &  23728 & 383578  & 668.46  & 295.51 & 3.2e+05 &  3.13 &  55341 \\
        long15   &  32767 & 753676  & 812.86  &  31.54 & 7.6e+05 &  3.97 &  51486 \\
        netlarge1&  45774 & 7238591 & 289.08  & 345.61 & 1.7e+06 & 44.95 &  89712 \\
        netlarge2&  39893 & 1158027 & 1043.61 &  41.29 & 1.2e+06 &  9.41 &  69478 \\
        netlarge3&  38502 & 4489009 & 990.27  & 117.03 & 3.9e+06 & 24.50 &  69326 \\
        netlarge6&   8000 & 15000000 & 774.88  &  60.26 & 1.5e+07 & 19.33 &  11736 \\
        square15 &  32760 & 753512  & 1191.03 &  39.45 & 7.6e+05 &  3.83 &  53391 \\
        wide15   &  32767 & 753676  &  807.88 &  33.29 & 7.6e+05 &  3.86 &  51486 \\
  \bottomrule
    \end{tabular}}}
  \label{tab: complete mcf benchmark results}
\end{table}%

\begin{table}[htbp]
  \centering
  \caption{Complete results of \Cref{fig: goto}, GOTO\_8}
  \setlength{\tabcolsep}{5mm}{
  \scalebox{0.8}{
    \begin{tabular}{lrrrrrrr}
    \toprule
    \multicolumn{1}{c}{\multirow{2}[3]{*}{prob}} & \multicolumn{1}{c}{\multirow{2}[3]{*}{node}} & \multicolumn{1}{c}{\multirow{2}[3]{*}{arc}} & \multicolumn{1}{c}{\multirow{2}[3]{*}{grbBarrier}} & \multicolumn{2}{c}{grbCrossover} & \multicolumn{2}{c}{CNET-grb} \\
\cmidrule{5-8}          &       &       &       & \multicolumn{1}{c}{time} & \multicolumn{1}{c}{iteration} & \multicolumn{1}{c}{time} & \multicolumn{1}{c}{iteration} \\
    \midrule
    8\_13a & 8192  & 65536 & 6.889 & 9.18  & 56895 & 1.13  & 15007 \\
    8\_13b & 8192  & 65536 & 7.65  & 9.62  & 53017 & 0.85  & 14448 \\
    8\_13c & 8192  & 65536 & 8.96  & 10.42  & 45673 & 0.54  & 11044 \\
    8\_13d & 8192  & 65536 & 6.7   & 9.53  & 53242 & 0.85  & 1409 \\
    8\_13e & 8192  & 65536 & 8.41  & 10.48  & 53550 & 0.93  & 14943 \\
    \midrule
    8\_14a & 16384 & 131072 & 19.22 & 26.17 & 108311 & 2.59  & 27530 \\
    8\_14b & 16384 & 131072 & 24.31 & 31.39 & 107992 & 2.66  & 29288 \\
    8\_14c & 16384 & 131072 & 15.36 & 21.84 & 102828 & 2.21  & 25294 \\
    8\_14d & 16384 & 131072 & 16.79 & 22.77 & 100342 & 2.47  & 25109 \\
    8\_14e & 16384 & 131072 & 21.82 & 28.15 & 102061 & 2.67  & 28538 \\
    \midrule
    8\_15a & 32768 & 262144 & 159.43 & 209.64 & 267968 & 4.07  & 60136 \\
    8\_15b & 32768 & 262144 & 147.3 & 185.93 & 262017 & 3.73  & 57446 \\
    8\_15c & 32768 & 262144 & 177.38 & 180.61 & 244905 & 2.48  & 44871 \\
    8\_15d & 32768 & 262144 & 150.96 & 152.57 & 202674 & 2.44  & 44258 \\
    8\_15e & 32768 & 262144 & 165.39 & 213.27 & 262002 & 6.44  & 60810 \\
    \midrule
    8\_16a & 65536 & 524288 & 130.68 & 256.75 & 441523 & 52.47 & 112601 \\
    8\_16b & 65536 & 524288 & 143.36 & 255.64 & 413308 & 47.15 & 106872 \\
    8\_16c & 65536 & 524288 & 121.97 & 261.29 & 422912 & 53.19 & 122950 \\
    8\_16d & 65536 & 524288 & 129.88 & 265.76 & 424800 & 53.4  & 123567 \\
    8\_16e & 65536 & 524288 & 97.37 & 300.54 & 519098 & 72.13 & 164305 \\
    \midrule
    8\_17a & 131072 & 1048576 & 318.76 & 773.08 & 883037 & 241.51 & 234596 \\
    8\_17b & 131072 & 1048576 & 324.16 & 872.61 & 872017 & 358.82 & 280343 \\
    8\_17c & 131072 & 1048576 & 262.28 & 796.8 & 865727 & 324.32 & 270919 \\
    8\_17d & 131072 & 1048576 & 248.78 & 892.76 & 916503 & 364.76 & 307621 \\
    8\_17e & 131072 & 1048576 & 317.87 & 1006.43 & 973954 & 446.39 & 339690 \\
    \bottomrule
    \end{tabular}}}%
\end{table}%

\begin{table}[htbp]
  \centering
  \caption{Complete results of \Cref{fig: goto}, GOTO\_sr}
  \setlength{\tabcolsep}{5mm}{
  \scalebox{0.8}{
    \begin{tabular}{lrrrrrrr}
    \toprule
    \multicolumn{1}{c}{\multirow{2}[3]{*}{prob}} & \multicolumn{1}{c}{\multirow{2}[3]{*}{node}} & \multicolumn{1}{c}{\multirow{2}[3]{*}{arc}} & \multicolumn{1}{c}{\multirow{2}[3]{*}{grbBarrier}} & \multicolumn{2}{c}{grbCrossover} & \multicolumn{2}{c}{CNET-grb} \\
\cmidrule{5-8}          &       &       &       & \multicolumn{1}{c}{time} & \multicolumn{1}{c}{iteration} & \multicolumn{1}{c}{time} & \multicolumn{1}{c}{iteration} \\
    \midrule
    sr\_12a & 4096  & 262144 & 9.11  & 4.13  & 61338 & 0.72  & 9257 \\
    sr\_12b & 4096  & 262144 & 11.36 & 8.45  & 74070 & 1.13  & 12128 \\
    sr\_12c & 4096  & 262144 & 9.92  & 5.91  & 73286 & 0.93  & 12464 \\
    sr\_12d & 4096  & 262144 & 8.93  & 10.95 & 69361 & 1.17  & 12790 \\
    sr\_12e & 4096  & 262144 & 8.81  & 6.98  & 65259 & 1.03  & 10744 \\
    \midrule
    sr\_13a & 8192  & 741455 & 24.64 & 94.02 & 283373 & 9.76  & 39252 \\
    sr\_13b & 8192  & 741455 & 22.01 & 39.54 & 280913 & 2.5   & 21259 \\
    sr\_13c & 8192  & 741455 & 28.57 & 49.29 & 256458 & 5.37  & 30001 \\
    sr\_13d & 8192  & 741455 & 29.86 & 91.38 & 277565 & 9.53  & 39677 \\
    sr\_13e & 8192  & 741455 & 29.34 & 93.98 & 274104 & 9.52  & 41648 \\
    \midrule
    sr\_14a & 16384 & 2097152 & 94.78 & 773.83 & 596187 & 41.47 & 93324 \\
    sr\_14b & 16384 & 2097152 & 127.57 & 694   & 595303 & 42.28 & 87569 \\
    sr\_14c & 16384 & 2097152 & 148.57 & 598.3 & 602842 & 41.21 & 85370 \\
    sr\_14d & 16384 & 2097152 & 190.37 & 127.52 & 522578 & 7.74  & 40602 \\
    sr\_14e & 16384 & 2097152 & 248.13 & 156.52 & 574786 & 13.7  & 56660 \\
    \midrule
    sr\_15a & 32768 & 5931642 & 778.51 & t     & -     & 203.16 & 193826 \\
    sr\_15b & 32768 & 5931642 & 831.5 & t     & -     & 176.6 & 152008 \\
    sr\_15c & 32768 & 5931642 & 1582.08 & 966.51 & 1862371 & 40.06 & 89170 \\
    sr\_15d & 32768 & 5931642 & 758.11 & t     & -     & 202.62 & 190106 \\
    sr\_15e & 32768 & 5931642 & 747.57 & t     & -     & 195.94 & 178075 \\
    \midrule
    sr\_16a & 65536 & 16777216 & 1715.79 & t     & -     & 1796  & 437039 \\
    sr\_16b & 65536 & 16777216 & 1605.52 & t     & -     & 1701.88 & 478398 \\
    sr\_16c & 65536 & 16777216 & 1677.75 & t     & -     & 1821.62 & 404442 \\
    sr\_16d & 65536 & 16777216 & 1376.29 & t     & -     & 1840.97 & 453485 \\
    sr\_16e & 65536 & 16777216 & 1310.3 & t     & -     & 1868  & 420691 \\
    \bottomrule
    \end{tabular}}}%
\end{table}%

\clearpage
\newpage

\section{Complete Numerical Results of Perturbation Crossover Method}\label{app: perturb experiment}

\paragraph{Problem set.} Our experiment includes all problems from Hans Mittelmann's benchmark, with the exclusion of those that encountered failures or timeouts during the barrier algorithm phase. 

Notably, \texttt{neos-5251015} and \texttt{physiciansched3-3} do not proceed as expected in \textsc{Gurobi}. For the former, \textsc{Gurobi} automatically switches to the homogeneous algorithm and for the latter, \textsc{Gurobi} yields a sub-optimal solution and terminates. Similarly, \texttt{ns1687037} fails in \textsc{Mosek}, which reports the solution is not found. \texttt{ns1688926} fails in both \textsc{Gurobi} and \textsc{Mosek}, with the former reporting a switch to the homogeneous algorithm and the latter reporting the solution is not found. Additionally, some instances like \texttt{thk\_48}, \texttt{thk\_63}, and \texttt{fhnw-binschedule1} are too large to complete within the one-hour barrier algorithm limit. 

The results for all the other problems are detailed in \Cref{tab: grb perturb complete} and \Cref{tab: msk perturb complete}.

Moreover, we note that among all the tested problems, \texttt{datt256\_lp} and \texttt{ex10} are feasibility problems. Consequently, we apply random perturbations exclusively to these two problems (see \Cref{fig:figureofalg}).

 Additionally, for some instances, the initial estimated optimal face might be infeasible. In those cases, we progressively decrease $\gamma$ (from 1e-3 to 1e-8, 1e-13, 1e-18, etc) to enlarge the estimated optimal face (see \Cref{fig:figureofalg}). This adjustment is necessary for 5 instances during our experiment: \texttt{cont11}, \texttt{dlr1}, \texttt{set-cover-model}, and \texttt{irish-electricity} (applied to both \textsc{Gurobi} and \textsc{Mosek}), and \texttt{Linf\_520c} (specific to \textsc{Mosek} only). For \texttt{irish-electricity}, a feasible candidate optimal face is achieved at $\gamma = $ 1e-13, whereas for the remaining 4 instances, it is $\gamma =$ 1e-8.

\begin{table}[htbp]
  \centering
  \caption{
  Complete result of the \textsc{Gurobi} part in \Cref{tab: perturb statistics}. The columns `$\text{P}_{\text{bfs}}$ Time' and '$\text{PD}_{\text{bfs}}$ Time' represent the running time needed to obtain an optimal primal BFS within tolerance $\expnumber{1}{-8}$ and the running time needed to obtain primal and dual optimal BFSs. }
    \begin{threeparttable}
    \setlength{\tabcolsep}{5mm}{
    \scalebox{0.7}{  
      \begin{tabular}{lrrrrrrr}
        \toprule
        \multicolumn{1}{c}{\multirow{2}{*}{Problem\tnote{1}}} & \multicolumn{1}{c}{\multirow{2}{*}{Barrier\tnote{2}}} & \multicolumn{1}{c}{\multirow{2}{*}{Crossover\tnote{3}}} & \multicolumn{4}{c}{Perturbation Crossover} \\
        \cmidrule(r){4-7}
        \multicolumn{1}{c}{} & \multicolumn{1}{c}{} & \multicolumn{1}{c}{} & \multicolumn{1}{c}{Perturb. Time\tnote{4}} & \multicolumn{1}{c}{Gap\tnote{5}} & \multicolumn{1}{c}{$\text{P}_{\text{bfs}}$ Time} & \multicolumn{1}{c}{$\text{PD}_{\text{bfs}}$ Time} \\
        \midrule
         s82                 &         143.38 &           881.19 &      0.53 &        2.2e-11 &    0.53 &  2329.26 \\
         datt256\_lp          &           1.65 &           415.94 &     18.19 &        1.2e-10 &   18.19 &    21.66 \\
         graph40-40\_lp       &           0.38 &           327.07 &     47.15 &        1.1e-13 &   47.15 &    47.34 \\
         set-cover-model     &         132.85 &           282.14 &      1.28 &        2.3e-10 &    1.28 &   103.79 \\
         woodlands09         &           6.59 &           196.39 &    157.94 &        1.9e-12 &  157.94 &   158.08 \\
         a2864               &           0.23 &           147.84 &    144.48 &        3.2e-13 &  144.48 &   229.58 \\
         nug08-3rd           &           0.87 &            55.7  &      0.01 &        3.6e-12 &    0.01 &   566.2  \\
         savsched1           &           7.17 &            54.14 &     39.8  &        2.4e-11 &   39.8  &    40.13 \\
         karted              &           6.55 &            46.72 &     55.6  &        1.7e-10 &   55.6  &    85.48 \\
         degme               &          22.05 &            42.77 &     46.8  &        4.6e-10 &   46.8  &    61.51 \\
         dlr1                &          89.66 &            33.73 &     73.12 &        3.4e-02   & 1511.18 &  1511.18\\
         neos3               &           1.83 &            23.74 &      2.92 &        3.2e-11 &    2.92 &  1555.26 \\
         s100                &          17.62 &            23.73 &      0.06 &        1.2e-12 &    0.06 &   116.6  \\
         rmine15             &          56.06 &            22.01 &      0.59 &        6.4e-10 &    0.59 &   187.88 \\
         tpl-tub-ws1617      &          55.08 &            12.76 &     12.07 &        3.7e-09 &   12.07 &    58.73 \\
         cont11              &           3.47 &            11.12 &     18.92 &        2.6e-07 &  841.23 &   841.23 \\
         cont1               &           2.39 &             9.05 &     10.57 &        2.5e-11 &   10.57 &    12.68 \\
         ex10                &           9.51 &             5.34 &     90.86 &        2.9e-11 &   90.86 &    90.9  \\
         ns1687037           &           6.91 &             5.31 &    178.69 &        3e-11   &  178.69 &   179.13 \\
         supportcase19       &          30.83 &             2.41 &      0.18 &        3.9e-09 &    0.18 &   237.49 \\
         rail02              &          41.33 &             1.95 &      0.3  &        1.5e-10 &    0.3  &  2969.26 \\
         rail4284            &          26.74 &             1.84 &      0.34 &        7.6e-10 &    0.34 &     1.33 \\
         Linf\_520c           &           9.86 &             1.8  &      6.82 &        7.2e-10 &    6.82 &     6.91 \\
         stormG2\_1000        &          88.8  &             1.66 &      3.05 &        1.9e-12 &    3.05 &     8.45 \\
         square41            &           1.24 &             1.52 &      0.05 &        9.5e-10 &    0.05 &    14.55 \\
         qap15               &           0.57 &             1.32 &      3.8  &        2.7e-09 &    3.8  &    22.76 \\
         scpm1\_lp            &          21.59 &             1.28 &      3.23 &        3.4e-10 &    3.23 &     3.42 \\
         pds-100             &          36.57 &             1.22 &      1.33 &        1.4e-11 &    1.33 &    15.6  \\
         stp3d               &           7.92 &             0.89 &      0.08 &        9.4e-12 &    0.08 &    81.51 \\
         fome13              &           1.1  &             0.89 &      0.16 &        1.9e-10 &    0.16 &    17.2  \\
         shs1023             &          28.42 &             0.67 &      0.19 &        1.4e-07 &   23.33 &    23.33 \\
         irish-electricity   &           4.5  &             0.51 &      4.68 &        6.7e-02 &   54.3  &    54.3  \\
         neos-5052403-cygnet &           4.78 &             0.44 &      1.06 &        2.2e-10 &    1.06 &    16.75 \\
         s250r10             &          12.81 &             0.3  &      0.07 &        1.5e-09 &    0.07 &     1.12 \\
         neos                &          15.06 &             0.22 &      0.28 &        2.4e-14 &    0.28 &     0.58 \\
         L1\_sixm250obs       &           2.14 &             0.02 &      0.05 &        6.6e-10 &    0.05 &     0.06 \\
        \bottomrule
        \end{tabular}
        }}
    \begin{tablenotes}
    \scriptsize
     \item[1] The problems are sorted by the running time of the crossover method of the solver.
     \item[2] \textsc{Gurobi}'s barrier method on the original problem without crossover stage;
     \item[3] \textsc{Gurobi}'s crossover stage on the original problem;
     \item[4] The time for solving the perturbed model via \textsc{Gurobi}
     \item[5] The relative gap between the primal (from the perturbed model) and the dual (from the barrier).
     \end{tablenotes}
  \end{threeparttable}
  \label{tab: grb perturb complete}
\end{table}

\begin{table}[htbp]
  \centering
  \caption{Complete result of the \textsc{Mosek} part in \Cref{tab: perturb statistics}. The columns `$\text{P}_{\text{bfs}}$ Time' and '$\text{PD}_{\text{bfs}}$ Time' represent the running time needed to obtain an optimal primal BFS within tolerance $\expnumber{1}{-8}$ and the running time needed to obtain primal and dual optimal BFSs.}
  \begin{threeparttable}
    \setlength{\tabcolsep}{5mm}
    \scalebox{0.7}{
       \begin{tabular}{lrrrrrrr}
        \toprule
        \multicolumn{1}{c}{\multirow{2}{*}{Problem\tnote{1}}} & \multicolumn{1}{c}{\multirow{2}{*}{Barrier\tnote{2}}} & \multicolumn{1}{c}{\multirow{2}{*}{Crossover\tnote{3}}} & \multicolumn{4}{c}{Perturbation Crossover} \\
        \cmidrule(r){4-7}
        \multicolumn{1}{c}{} & \multicolumn{1}{c}{} & \multicolumn{1}{c}{} & \multicolumn{1}{c}{Perturb. Time\tnote{4}} & \multicolumn{1}{c}{Gap\tnote{5}} & \multicolumn{1}{c}{$\text{P}_{\text{bfs}}$ Time} & \multicolumn{1}{c}{$\text{PD}_{\text{bfs}}$ Time} \\
        \midrule
         cont11              &           2.64 &          2950.6  &   1070.38 &        4.5e-01    & t\tnote{6} &  t \\
         a2864               &          61.5  &           555.86 &   t  &      nan       &  t    &   t    \\
         tpl-tub-ws1617      &          36.24 &           315.11 &      5.42 &        7.3e-07 & 1300.5  &  1300.5  \\
         datt256\_lp          &           1.38 &           205.85 &      3.14 &        1.8e-12 &    3.14 &     4.56 \\
         woodlands09         &          11.41 &           143.76 &    159.64 &        4.2e-09 &  159.64 &  t \\
         cont1               &           3.45 &            91.08 &   2534.38 &      1.79e-12   & 2534.38 &  2634.29 \\
         nug08-3rd           &          14.12 &            69.79 &      0.01 &        4.9e-12 &    0.01 &  3604.22 \\
         savsched1           &           8.45 &            60.33 &   1578.21 &        1.4e-12 & 1578.21 &  2691.28 \\
         graph40-40\_lp       &           3.84 &            42.82 &     17.59 &        3.5e-13 &   17.59 &    18.58 \\
         degme               &         140.89 &            38.02 &     77.35 &        1.5e-09 &   77.35 &   224.08 \\
         dlr1                &         176.46 &            37.09 &    249.78 &        1.3e-01    & t &  t \\
         karted              &          96.22 &            27.58 &     35.24 &        3.9e-09 &   35.24 &   134.57 \\
         set-cover-model     &         176.38 &            12.76 &      1.14 &        1.1e-10 &    1.14 &    15.07 \\
         Linf\_520c           &           9.27 &             7.97 &    160.39 &        9.4e-01    & 1940.92 &  1940.92 \\
         neos3               &           5.26 &             6.87 &     85.89 &        1.3e-08 &  494.79 &   494.79 \\
         rmine15             &          47.85 &             5.93 &      3.81 &        2.4e-10 &    3.81 &  1213.63 \\
         ex10                &           7.12 &             5.41 &     13.24 &        2.8e-08 &  t    &   t    \\
         fome13              &           4.54 &             3.95 &      0.34 &        9.7e-09 &    0.34 &   213.4  \\
         s82                 &          45.44 &             2.38 &      0.29 &        1.3e-10 &    0.29 &  1225.62 \\
         stormG2\_1000        &          26.7  &             2.3  &      1.87 &        4.9e-08 &   26.86 &    26.86 \\
         supportcase19       &           9.42 &             1.96 &      0.07 &        1.5e-11 &    0.07 &  3343.68 \\
         pds-100             &          33.11 &             1.32 &      0.88 &        2e-10   &    0.88 &   114.26 \\
         scpm1\_lp            &          20.09 &             1.3  &      3.66 &        6.8e-11 &    3.66 &     6.06 \\
         qap15               &           4.58 &             1.12 &      1.63 &        9.4e-14 &    1.63 &     6.28 \\
         stp3d               &           6.85 &             1.06 &      0.08 &        8.2e-12 &    0.08 &   639.58 \\
         square41            &           4.21 &             1.04 &      0.1  &        5.5e-15 &    0.1  &    48.92 \\
         rail4284            &          10.93 &             0.79 &      0.72 &        6.7e-14 &    0.72 &     4.57 \\
         rail02              &          46.01 &             0.61 &      0.47 &        1.5e-07 &  827.25 &   827.25 \\
         physiciansched3-3   &          22.27 &             0.58 &      5.01 &        4.1e-11 &    5.01 &    91.31 \\
         neos-5052403-cygnet &           3.68 &             0.53 &      1.11 &        4.4e-12 &    1.11 &     6.24 \\
         neos                &           9.16 &             0.48 &      0.47 &        2.9e-09 &    0.47 &     5.79 \\
         shs1023             &          23.63 &             0.41 &      0.14 &        1.9e-08 &   30.49 &    30.49 \\
         irish-electricity   &           7.3  &             0.38 &     11.17 &        8.8e-02   &  243.57 &   243.57 \\
         s100                &           4.9  &             0.12 &      0.02 &        3.6e-13 &    0.02 &    39.61 \\
         s250r10             &           3.47 &             0.09 &      0.05 &        2.2e-12 &    0.05 &    14.51 \\
         L1\_sixm250obs       &           3.17 &             0.02 &      0.03 &        1e-10   &    0.03 &     0.07 \\
        \bottomrule
      \end{tabular}
    }
    \begin{tablenotes}
     \scriptsize
     \item[1] The problems are sorted by the running time of the crossover method of the solver.
     \item[2] \textsc{Mosek}'s barrier method on the original problem without crossover stage;
     \item[3] \textsc{Mosek}'s crossover stage on the original problem;
     \item[4] The time for solving the perturbed model via \textsc{Mosek}
     \item[5] The relative gap between the primal (from the perturbed model) and the dual (from the barrier).
     \item[6] Exceeds time limit: 3600 seconds.
    \end{tablenotes}
  \end{threeparttable}
  \label{tab: msk perturb complete}
\end{table}

\clearpage

\section{A Comparison between Different Perturbation Methods}\label{app: different perturb}

In this section, we compare different perturbation crossover methods, including an apriori perturbation method, \citet{mehrotra1991finding}'s perturbation method, and our method. An apriori perturbation aims at perturbing the cost vector to avoid multiple optimal solutions, before running the barrier algorithm. \citet{mehrotra1991finding}'s method uses an interior-point solution to estimate a better perturbation, while ours added more techniques to control the magnitude of the perturbation and restricted the problem to an estimated optimal face. 

\begin{table}[h!]
 \centering
 \caption{Differences among perturbation methods.}
 \vspace{+3pt}
    \setlength{\tabcolsep}{0.1mm}{
    \scalebox{0.8}{
        \begin{tabular}{cccc}
            \toprule
                & \begin{tabular}{c} Use info from an \\ interior-point solution $x$\end{tabular}  & Perturbation strategy & Estimate the optimal face  \\
            \midrule
             \textsc{Apriori perturb} & No &  Random  &  No \\
             \textsc{Mehrotra's method} & Yes & Random + info of $x$ & No \\
             \textsc{Our method} & Yes & Random + info of $x$ and problem data & Yes \\
             \bottomrule
        \end{tabular}
    }
    }
    \label{tab: compare perturb methods (def)}
\end{table}
 
We test these methods in the same problem set as we did in \Cref{subsec: perturbation experiments}. Mehrotra's method can be referred to as (2.6) in \citet{mehrotra1991finding} and we set the $\epsilon$ to be 1. 
In  Mehrotra's method,
we additionally set the perturbation corresponding to tiny components of the given solution as zero to avoid ill-condition cases. 
\Cref{tab: compare perturb methods} shows the significant advantage of our perturbation crossover method. 

\begin{table}[h!]
  \centering
    \caption{Performance statistics of different perturbation crossover methods: an apriori random perturbation, Mehrotra's method \citep{mehrotra1991finding}, and our method using \textsc{Gurobi}. The goals are to find a primal BFS within precision $\expnumber{1}{-8}$ for all three methods. The numbers of total tested problems are 36 for all three methods. The time limit is set to 500 seconds for this experiment.}
    \vspace{5pt}
  \begin{threeparttable}
    \setlength{\tabcolsep}{3mm}{
    \scalebox{0.9}{
      \begin{tabular}{lcccc}
            \toprule
                & Faster & Slower & Timeout  & GeoAvg Time \\
            \midrule
                \textsc{Apriori perturbation}  & 6  & 30  &  5 & 53.41    \\
                \textsc{Mehrotra's method}  & 2  & 34  & 4  & 40.75    \\
                \textsc{Our method}         & 19 & 17  & 2 &  3.38    \\
            \bottomrule
        \end{tabular}}}
    \vspace{+0.25cm}
  \end{threeparttable}
  \label{tab: compare perturb methods}
\end{table}

\end{APPENDICES}

\end{document}